\documentclass{article}

\usepackage{amsfonts,amsmath,amssymb}
\usepackage{hyperref}
\usepackage{color}
\usepackage{graphicx}
\usepackage[section]{placeins}
\usepackage[utf8]{inputenc}
\usepackage{tikz}

\definecolor{darkblue}{rgb}{0,0.08,0.45}
\hypersetup{colorlinks,linkcolor=darkblue,citecolor=darkblue,filecolor=darkblue,urlcolor=darkblue}
\hypersetup{bookmarksopen,bookmarksnumbered}

\numberwithin{equation}{section}

\renewcommand{\d}{\ensuremath{\,\mathrm{d}}}

\newcommand{\normal}{\ensuremath{\hat{\mathbf{n}}}}
\newcommand{\utot}{\ensuremath{p_{\mathrm{tot}}}}
\newcommand{\usca}{\ensuremath{p_{\mathrm{sca}}}}
\newcommand{\uinc}{\ensuremath{p_{\mathrm{inc}}}}
\newcommand{\SLP}{\ensuremath{\mathcal{V}}}
\newcommand{\DLP}{\ensuremath{\mathcal{K}}}

\newcommand{\SL}{\ensuremath{V}}
\newcommand{\DL}{\ensuremath{K}}
\newcommand{\AD}{\ensuremath{T}}
\newcommand{\HS}{\ensuremath{D}}
\newcommand{\ID}{\ensuremath{I}}
\newcommand{\CA}{\ensuremath{A}}

\newcommand{\NtDmi}{\ensuremath{\Lambda_{\mathrm{NtD},m}^-}}
\newcommand{\DtNmi}{\ensuremath{\Lambda_{\mathrm{DtN},m}^-}}

\newcommand{\OsrcNtDmi}{\ensuremath{L_{\mathrm{NtD},m}^-}}
\newcommand{\OsrcDtNmi}{\ensuremath{L_{\mathrm{DtN},m}^-}}

\newcommand{\traceDe}{\ensuremath{\gamma_D^+}}

\newcommand{\traceNe}{\ensuremath{\gamma_N^+}}

\newcommand{\traceDm}{\ensuremath{\gamma_{D,m}}}

\newcommand{\traceDme}{\ensuremath{\gamma_{D,m}^+}}

\newcommand{\traceDmi}{\ensuremath{\gamma_{D,m}^-}}

\newcommand{\traceNm}{\ensuremath{\gamma_{N,m}}}

\newcommand{\traceNme}{\ensuremath{\gamma_{N,m}^+}}

\newcommand{\traceNmi}{\ensuremath{\gamma_{N,m}^-}}

\title{Benchmarking preconditioned boundary integral formulations for acoustics\footnote{This is the peer reviewed version of the following article: ``Benchmarking preconditioned boundary integral formulations for acoustics,'' the \emph{International Journal for Numerical Methods in Engineering}, volume~122, issue~20, pages~5873--5897 (2021), which has been published in final form at \url{https://doi.org/10.1002/nme.6777}.\newline This article may be used for non-commercial purposes in accordance with Wiley Terms and Conditions for Use of Self-Archived Versions. This article may not be enhanced, enriched or otherwise transformed into a derivative work, without express permission from Wiley or by statutory rights under applicable legislation. Copyright notices must not be removed, obscured or modified. The article must be linked to Wiley’s version of record on Wiley Online Library and any embedding, framing or otherwise making available the article or pages thereof by third parties from platforms, services and websites other than Wiley Online Library must be prohibited.}}
\author{Elwin van 't Wout\thanks{Institute for Mathematical and Computational Engineering, School of Engineering and Faculty of Mathematics, Pontificia Universidad Católica de Chile, Santiago, Chile. Contact: e.wout@uc.cl} \and 
Seyyed R.~Haqshenas\thanks{Department of Mechanical Engineering, University College London, London, United Kingdom.} \and
Pierre Gélat\footnotemark[2] \and
Timo Betcke\thanks{Department of Mathematics, University College London, London, United Kingdom.} \and
Nader Saffari\footnotemark[2]}
\date{October 31, 2022}

\begin{document}

\maketitle

\begin{abstract}
	The boundary element method (BEM) is an efficient numerical method for simulating harmonic wave propagation. It uses boundary integral formulations of the Helmholtz equation at the interfaces of piecewise homogeneous domains. The discretisation of its weak formulation leads to a dense system of linear equations, which is typically solved with an iterative linear method such as GMRES. The application of BEM to simulating wave propagation through large-scale geometries is only feasible when compression and preconditioning techniques reduce the computational footprint. Furthermore, many different boundary integral equations exist that solve the same boundary value problem. The choice of preconditioner and boundary integral formulation is often optimised for a specific configuration, depending on the geometry, material characteristics, and driving frequency. On the one hand, the design flexibility for the BEM can lead to fast and accurate schemes. On the other hand, efficient and robust algorithms are difficult to achieve without expert knowledge of the BEM intricacies. This study surveys the design of boundary integral formulations for acoustics and their acceleration with operator preconditioners. Extensive benchmarks provide valuable information on the computational characteristics of several hundred different models for multiple reflection and transmission of acoustic waves.
\end{abstract}

\section{Introduction}

The Helmholtz equation for harmonic wave propagation is a widely used model for many acoustic phenomena, such as room acoustics, sonar, and biomedical ultrasound, among others~\cite{lahaye2017modern}. The boundary element method (BEM) is one of the most efficient numerical methods to solve Helmholtz transmission problems and is based on boundary integral formulations that rewrite the volumetric partial differential equations into a representation of the acoustic fields in terms of surface potentials at the material interfaces~\cite{nedelec2001acoustic, steinbach2008numerical, hsiao2008boundary, sauter2010boundary}. Many different boundary integral formulations can model precisely the same physical problem. This design flexibility allows for the development of specialised formulations but causes complications for many practitioners who are obliged to pick a formulation from decades of scientific literature or go through the intrinsic design process themselves. Furthermore, modern preconditioning techniques that considerably improve the BEM's computational efficiency have yet to be applied to many. Here, we will present five different design strategies, develop efficient preconditioners, and compare the computational characteristics of hundreds of preconditioned boundary integral formulations through extensive benchmarking of acoustic transmission at multiple penetrable domains.

The BEM has unique advantages over volumetric methodologies such as finite element and finite difference methods. Firstly, unbounded exterior domains are naturally handled since the representation formulas automatically satisfy the radiation conditions, thus avoiding artificial boundary conditions to truncate the computational domain. Secondly, the number of degrees of freedom scales quadratically with respect to the frequency. Thirdly, the fast multipole method~\cite{greengard1987fast} and hierarchical matrix compression~\cite{hackbusch2015hierarchical} perform dense matrix arithmetic in almost linear scaling. Fourthly, Green's functions are explicitly used and numerical dispersion or dissipation is expected to be limited~\cite{baydoun2018quantification}. Finally, open-source software provides high-level programming platforms~\cite{smigaj2015solving}. On the downside, the BEM is limited to problem settings for which Green's functions are available. For this reason, the geometry needs to consist of piecewise homogeneous materials. Considering the BEM's advantages and limitations, it is the preferred methodology to simulate many wave propagation problems with applications in acoustics, electromagnetics and elastodynamics~\cite{chew2001fast, marburg2018boundary, kirkup2019boundary}.

The BEM reformulates the Helmholtz equation into a boundary integral equation before the discretisation process. In contrast, volumetric methods discretise the Helmholtz equation directly. Since the boundary integral formulation uses potential theory, there is great flexibility in defining the fields' representation in terms of surface potentials. The many design strategies that are availabe lead to an infinite number of boundary integral formulations for the same acoustic transmission problem. An abundance of different formulations have been introduced in the scientific literature in the last decades: first for mathematical analysis of rigid scatterers~\cite{maue1949formulierung} and quickly extended to transmission into penetrable domains~\cite{muller1957grundprobleme, mitzner1966acoustic, mautz1977electromagnetic, costabel1985direct}. Most of the boundary integral formulations are presented with different notational frameworks, are designed through different processes, and are often dedicated to specific application areas, frequency ranges, material types or discretisation techniques. Furthermore, techniques such as robust singular integration, fast multipole methods, hierarchical matrix compression, parallel computing and preconditioning have improved computational efficiency tremendously over the last decades. Hence, formulations that were inefficient when introduced in literature might have become competitive with modern-day algorithms.

While the high level of design flexibility for the BEM is beneficial to the expert who can design efficient algorithms for a specific purpose, it is a burden to many practitioners who need to find the correct mathematical framework and computational configuration for their simulation settings. This study summarises five families of boundary integral formulations for acoustic transmission through multiple domains. Three different families (single-trace, multiple-traces and auxiliary field formulations) use direct representation formulas for the acoustic field, each for a different set of surface potentials. The other two families (single potential and mixed potential formulations) use indirect representations for either all fields or only the exterior fields, respectively. For the first time, operator preconditioning will be applied to all of these formulations. The main novelty of this study is the extensive benchmarking. This will provide insights into the computational performance of the preconditioned boundary integral formulations and their multifaceted dependencies on frequency range, material type and geometry. Different models will be compared in terms of calculation time, accuracy and convergence at canonical test cases and large-scale simulations.

The Helmholtz transmission problem will be detailed in Section~\ref{sec:formulation} along with the boundary integral operators. Section~\ref{sec:bie} then surveys most of the boundary integral formulations from the literature, and operator preconditioning is discussed in Section~\ref{sec:preconditioning}. The computational results from extensive benchmarking are presented in Section~\ref{sec:results}, followed by a discussion and conclusions on the study.

\section{Formulation}
\label{sec:formulation}

\subsection{Helmholtz equation}

\begin{figure}[!t]
	\centering
	\begin{tikzpicture}
		\draw[thin](-4.8,-1)--(-4.8,1);
		\draw[thin](-4.6,-1)--(-4.6,1);
		\draw[thin](-4.4,-1)--(-4.4,1);
		\draw[thick, ->](-5.2,0)--(-4,0);
		\node at (-3.9,-0.3) {$\uinc$};
		
		\node at (-2.5,-1) {$\Omega_0$};
		
		\draw[semithick](-2,1) circle (0.8) node {$\Omega_1$};
		\node at (-3,0.6) {$\Gamma_1$};
		\draw[thick, ->](-1.4343,1.5657)--(-1.1343,1.8657) node[above] {$\normal_1$};
		
		\draw[semithick](0,-1) ellipse (1.2 and 0.8) node {$\Omega_2$};
		\node at (-0.9,-0.1) {$\Gamma_2$};
		\draw[thick, ->](0,-0.2)--(0,0.2) node[right] {$\normal_2$};
		
		\draw[semithick](2,1) ellipse (0.8 and 1.0) node {$\Omega_3$};
		\node at (1.7,-0.2) {$\Gamma_3$};
		\draw[thick, ->](1.2,1)--(0.7,1) node[above] {$\normal_3$};
	\end{tikzpicture}
	\caption{A sketch of the geometry setting for the wave propagation model.}
	\label{fig:geometry}
\end{figure}
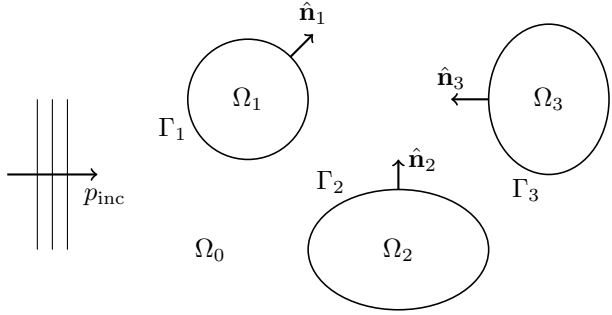

The Helmholtz equation is the standard model for the propagation of harmonic acoustic waves in materials with a linear response in the frequency domain. The geometry consists of a collection of objects embedded in free space, as depicted in Figure~\ref{fig:geometry}. Let us denote the exterior unbounded domain by~$\Omega_0 \subset \mathbb{R}^3$, and the objects by $\Omega_1, \Omega_2, \dots, \Omega_\ell$ for $\ell \geq 1$ a constant. All objects are assumed to be disjoint, bounded, and with a homogeneous interior, thus excluding junctions of interfaces. The wavenumber in each domain is denoted by $k_0, k_1, k_2, \dots, k_\ell$, respectively and each domain is equiped with a material constant $\sigma_n$, $n=0,1,\dots,\ell$ that typically depends on the mass density or acoustic impedance. Let us denote the boundaries of the objects by $\Gamma_1, \Gamma_2, \dots, \Gamma_\ell$ and assume they are smooth and can be equiped with unit normal vectors $\normal_1, \normal_2, \dots, \normal_\ell$ all pointing towards the exterior domain.

After extracting the time-dependent $e^{-\imath\omega t}$ term with $\omega$ the angular frequency, the harmonic acoustic pressure field is denoted by $\utot$ and is a complex-valued function. In the exterior, the pressure can be decomposed into the known incident and unknown scattered field as $\utot = \uinc + \usca$. The incident field can be any acoustic field that satisfies the Helmholtz equation with the wavenumber given by the exterior region and will be chosen to be a plane wave field in this study. The equations of motion for scalar harmonic wave propagation are given by the Helmholtz system as
\begin{equation} \label{eq:helmholtz}
	\begin{cases}
		\Delta \utot + k_m^2 \utot = 0, & \text{in } \Omega_m \text{ for } m=0,1,2,\dots,\ell; \\
		\traceDme \utot = \traceDmi \utot, & \text{at } \Gamma_m \text{ for } m=1,2,\dots,\ell; \\
		\sigma_0 \traceNme \utot = \sigma_m \traceNmi \utot, & \text{at } \Gamma_m \text{ for } m=1,2,\dots,\ell; \\
		\lim_{\mathbf{r} \to \infty} |\mathbf{r}| (\partial_{|\mathbf{r}|} \usca(\mathbf{r}) - ik_0 \usca(\mathbf{r}) = 0
	\end{cases}
\end{equation}
where continuity of the fields across the interfaces is assumed. The last equation is the Sommerfeld radiation condition which states that the scattered field radiates towards infinity, where $\mathbf{r}$ denotes the position.
The traces of the fields at the material interfaces are defined as
\begin{align}
	\traceDmi f(\mathbf{x}) &= \lim_{\Omega_m \ni \mathbf{y} \to \mathbf{x}} f(\mathbf{y}), \\
	\traceNmi f(\mathbf{x}) &= \lim_{\Omega_m \ni \mathbf{y} \to \mathbf{x}} \nabla f(\mathbf{y}) \cdot \normal_m(\mathbf{x}), \\
	\traceDme f(\mathbf{x}) &= \lim_{\Omega_0 \ni \mathbf{y} \to \mathbf{x}} f(\mathbf{y}), \\
	\traceNme f(\mathbf{x}) &= \lim_{\Omega_0 \ni \mathbf{y} \to \mathbf{x}} \nabla f(\mathbf{y}) \cdot \normal_m(\mathbf{x})
\end{align}
for $\mathbf{x} \in \Gamma_m$ and $m=1,2,\dots,\ell$, where the indices $\gamma_m^\pm$ indicate traces from the exterior or interior of subdomain~$m$, respectively. The traces~$\gamma_{D,N}$ are called the Dirichlet and Neumann traces, respectively, and are related to the acoustic pressure and normal particle velocity at the interface.

\subsection{Preliminaries for boundary integral formulations}

This section summarises the definitions and properties of boundary integral operators that will be used for the BEM. Proofs, details and more information can be found in standard literature~\cite{nedelec2001acoustic, steinbach2008numerical, hsiao2008boundary, sauter2010boundary}.

\subsubsection{Boundary integral operators}

The single-layer and double-layer potential integral operators that map from a surface potential at interface~$\Gamma_m$ towards the subdomain~$\Omega_j$ are given by
\begin{align}
	[\SLP_{j,m} \psi](\mathbf{x}) &= \iint_{\Gamma_m} G_j(\mathbf{x},\mathbf{y}) \psi(\mathbf{y}) \d\mathbf{y}, \\
	[\DLP_{j,m} \phi](\mathbf{x}) &= \iint_{\Gamma_m} \frac{\partial G_j(\mathbf{x},\mathbf{y})}{\partial \hat{\mathbf{n}}(\mathbf{y})} \phi(\mathbf{y}) \d\mathbf{y},
\end{align}
for $\mathbf{x} \in \Omega_j$.
Here, $G_j$ denotes the Green's function with the wavenumber of the respective region, that is,
\begin{equation}
	G_j(\mathbf{x},\mathbf{y}) = \frac{e^{\imath k_j |\mathbf{x}-\mathbf{y}|}}{4\pi|\mathbf{x}-\mathbf{y}|} \quad \text{for } \mathbf{x},\mathbf{y} \in \Omega_j \text{ and } \mathbf{x} \ne \mathbf{y}
\end{equation}
for $j=0,1,2,\dots,\ell$ where $\iota$ denotes the complex unit.
The boundary integral operators that map from one interface~$\Gamma_n$ to another or the same interface~$\Gamma_m$ are given by
\begin{align}
	[\SL_{j,mn} \psi](\mathbf{x}) &= \iint_{\Gamma_n} G_j (\mathbf{x},\mathbf{y}) \psi(\mathbf{y}) \d\mathbf{y} \\
	[\DL_{j,mn} \phi](\mathbf{x}) &= \iint_{\Gamma_n} \frac{\partial}{\partial \normal_n(\mathbf{y})} G_j (\mathbf{x},\mathbf{y}) \phi(\mathbf{y}) \d\mathbf{y}, \\
	[\AD_{j,mn} \psi](\mathbf{x}) &= \frac{\partial}{\partial \normal_m(\mathbf{x})} \iint_{\Gamma_n} G_j (\mathbf{x},\mathbf{y}) \psi(\mathbf{y}) \d\mathbf{y}, \\
	[\HS_{j,mn} \phi](\mathbf{x}) &= -\frac{\partial}{\partial \normal_m(\mathbf{x})} \iint_{\Gamma_n} \frac{\partial}{\partial \normal_n(\mathbf{y})} G_j (\mathbf{x},\mathbf{y}) \phi(\mathbf{y}) \d\mathbf{y},
\end{align}
for $\mathbf{x} \in \Gamma_m$, which are called the single-layer, double-layer, adjoint double-layer and hypersingular boundary integral operators, respectively. A single subscript~$m$ will be used for interior operators, that is, $\SL_m = \SL_{m,mm}$ and similar for the other operators.
The identity operator acting on a surface potential at interface~$\Gamma_m$ is denoted by $\ID_m$ and
\begin{equation}
	\bar\ID_m = \begin{bmatrix} \ID_m & 0 \\ 0 & \ID_m \end{bmatrix}
\end{equation}
denotes the identity operator acting on a pair of surface potentials at interface~$\Gamma_m$.
Furthermore,
\begin{equation}
	\CA_{j,mn} = \begin{bmatrix} -\DL_{j,mn} & \SL_{j,mn} \\ \HS_{j,mn} & \AD_{j,mn} \end{bmatrix}
\end{equation}
denotes the Calderón boundary integral operator.

\subsubsection{Calderón identities}
\label{sec:calderonidentities}

The Calderón operator satisfies the projection property
\begin{equation}
	A_{j,mm}^2 = \frac14\bar\ID_m.
\end{equation}
Hence,
\begin{align}
	\SL_{j,mm} \HS_{j,mm} &= \frac14 \ID_m - \DL_{j,mm}^2, \label{eq:slhs} \\
	\HS_{j,mm} \SL_{j,mm} &= \frac14 \ID_m - \AD_{j,mm}^2 \label{eq:hssl}
\end{align}
which are called Calderón identities.

\subsubsection{The Neumann-to-Dirichlet map}
\label{sec:ntdmap}

The interior Neumann-to-Dirichlet (NtD) and Dirichlet-to-Neumann (DtN) maps are implicitly defined as
\begin{align}
	\traceDmi \utot &= \NtDmi \traceNmi \utot, \label{eq:ntd} \\
	\traceNmi \utot &= \DtNmi \traceDmi \utot \label{eq:dtn}
\end{align}
which are also known as the Poincaré-Steklov and Steklov-Poincaré operators, respectively.
These operators satisfy
\begin{align}
	\NtDmi &= \left(\frac12 \ID_m + \DL_m\right)^{-1} \SL_m = \HS_m^{-1} \left(\frac12 \ID_m - \AD_m\right), \label{eq:ntd:operator} \\
	\DtNmi &= \left(\frac12 \ID_m - \AD_m\right)^{-1} \HS_m = \SL_m^{-1} \left(\frac12 \ID_m + \DL_m\right). \label{eq:dtn:operator}
\end{align}
and, by symmetry of the operators, one has
\begin{align}
	\NtDmi \HS_m &= \frac12 \ID_m - \DL_m, \label{eq:ntd:prec} \\
	\DtNmi \SL_m &= \frac12 \ID_m + \AD_m. \label{eq:dtn:prec}
\end{align}
The exterior NtD and DtN maps have similar expressions in the case of a single object but the extensions to multiple reflection requires a global system~\cite{acosta2015surface, alzubaidi2016formulation}.

\section{Boundary integral formulations}
\label{sec:bie}

Generally speaking, boundary integral formulations can be classified into direct and indirect formulations. Whereas direct formulations use a field representation given by both the single-layer and double-layer potentials operators acting on traces of the field, the indirect formulations use a field respresentation in terms of arbitrary surface potentials.

\subsection{Direct boundary integral formulations}
\label{sec:formulation:direct}

The direct representation formula for the acoustic field is given by
\begin{align}
	\usca &= - \sum_{n=1}^\ell \left(\SLP_{0,n} (\traceNme \utot) - \DLP_{0,n} (\traceDme \utot)\right) & \text{in } \Omega_0,
	\label{eq:representation:exterior:direct} \\
	\utot &= \SLP_m (\traceNmi \utot) - \DLP_m (\traceDmi \utot) & \text{in } \Omega_m
	\label{eq:representation:interior:direct}
\end{align}
for $m=1,2,\dots,\ell$.
Different procedures exist to use the interface conditions for coupling the exterior and interior representations. The \emph{single-trace formulations} (STF) use a single pair of two field traces as unknown surface potentials, \emph{multiple-traces formulations} (MTF) use all four traces of the fields, and \emph{auxiliary field formulations} (AFF) reduce the formulation to a single surface potential at each interface.

\subsubsection{Single-trace formulations}
\label{sec:formulation:singletrace}

The principle behind single-trace formulations is to use the transmission conditions to eliminate half of the unknown potentials by defining a single set of Dirichlet and Neumann traces as
\begin{align}
	\phi_m &= \traceDme \utot = \traceDmi \utot, \\
	\psi_m &= \traceNme \utot = \frac{\sigma_m}{\sigma_0} \traceNmi \utot
\end{align}
for $m=1,2,\dots,\ell$, where the impedance ratio could have been defined at the exterior as well.
Then, the traces of the representation formulas~\eqref{eq:representation:exterior:direct}--\eqref{eq:representation:interior:direct} are given by the Calderón equations
\begin{align}
	\left(\frac12\bar\ID_m + \CA_{0,mm} \right)
	\begin{bmatrix} \phi_m \\ \psi_m \end{bmatrix} + \sum_{n=1, n \ne m}^\ell \CA_{0,mn}
	\begin{bmatrix} \phi_n \\ \psi_n \end{bmatrix}
	&= \begin{bmatrix} \traceDme \uinc \\ \traceNme \uinc \end{bmatrix},
	\label{eq:calderon:exterior:stf} \\
	\left(\frac12\bar\ID_m - \widehat{\CA}_m \right)
	\begin{bmatrix} \phi_m \\ \psi_m \end{bmatrix}
	&= \begin{bmatrix} 0 \\ 0 \end{bmatrix}
	\label{eq:calderon:interior:stf}
\end{align}
for $m=1,2,\dots,\ell$, where
\begin{equation}
	\label{eq:calderon:scaled}
	\widehat{\CA}_m = \begin{bmatrix}
		-\DL_m & \frac{\sigma_0}{\sigma_m} \SL_m \\
		\frac{\sigma_m}{\sigma_0} \HS_m & \AD_m
	\end{bmatrix}
\end{equation}
is a scaled interior Calderón matrix.
This is a linear system of $4\ell$ equations for $2\ell$ unknown potentials and the design of single-trace formulations follow different approaches to reduce the dimensionality.

\paragraph{Dirichlet formulation}
Selecting the Dirichlet traces of the representation formulas results in
\begin{equation}
	\label{eq:formulation:dirichlet}
	\begin{cases}
		\frac12 \phi_m + \sum_{n=1}^\ell \left(-\DL_{0,mn} \phi_n + \SL_{0,mn} \psi_n\right) = \traceDme \uinc, \\
		\frac12 \phi_m + \DL_{m,mm} \phi_m - \frac{\sigma_0}{\sigma_m} \SL_{m,mm} \psi_m = 0
	\end{cases}
\end{equation}
for $m=1,2,\dots,\ell$ which is called the \emph{Dirichlet formulation}.

\paragraph{Neumann formulation}
Selecting the Neumann traces of the representation formulas results in
\begin{equation}
	\label{eq:formulation:neumann}
	\begin{cases}
		\frac12 \psi_m + \sum_{n=1}^\ell \left(\HS_{0,mn} \phi_n + \AD_{0,mn} \psi_n\right) = \traceNme \uinc, \\
		\frac12 \psi_m - \frac{\sigma_m}{\sigma_0} \HS_{m,mm} \phi_m - \AD_{m,mm} \psi_m = 0
	\end{cases}
\end{equation}
for $m=1,2,\dots,\ell$ which is called the \emph{Neumann formulation}.

\paragraph{PMCHWT formulation}
Taking the difference of the exterior and interior traces of the representation formulas results in
\begin{equation}
	\label{eq:formulation:pmchwt}
	\widehat\CA_m \begin{bmatrix} \phi_m \\ \psi_m \end{bmatrix} + \sum_{n=1}^\ell \CA_{0,mn}
	\begin{bmatrix} \phi_n \\ \psi_n \end{bmatrix}
	= \begin{bmatrix} \traceDme \uinc \\ \traceNme \uinc \end{bmatrix}
\end{equation}
for $m=1,2,\dots,\ell$ which is called the \emph{PMCHWT formulation} (Poggio-Miller-Chang-Harrington-Wu-Tsai)~\cite{poggio1973integral, chang1974surface, wu1977scattering-bor}.

\paragraph{Müller formulation}
Taking the sum of the exterior and interior traces of the representation formulas results in
\begin{equation}
	\label{eq:formulation:muller}
	\begin{bmatrix} \phi_m \\ \psi_m \end{bmatrix} - \widehat\CA_m \begin{bmatrix} \phi_m \\ \psi_m \end{bmatrix} + \sum_{n=1}^\ell \CA_{0,mn}
	\begin{bmatrix} \phi_n \\ \psi_n \end{bmatrix}
	= \begin{bmatrix} \traceDme \uinc \\ \traceNme \uinc \end{bmatrix}
\end{equation}
for $m=1,2,\dots,\ell$ which is called the \emph{Müller formulation}~\cite{muller1957grundprobleme}.

\paragraph{Combined formulations}
In general, arbitrary linear combinations of the traces of the representation formulas can be taken~\cite{mitzner1966acoustic, mautz1977electromagnetic}. For constants $\eta_m^\pm$ and $\nu_m^\pm$ one can distinguish the following formulations.
The \emph{combined trace formulation} is given by
\begin{equation}
	\label{eq:formulation:combined:trace}
	\begin{cases}
		\eta_m^+ \left( \frac12 \phi_m + \sum_{n=1}^\ell \left(-\DL_{0,mn} \phi_n + \SL_{0,mn} \psi_n\right) \right) \\
		\quad + \eta_m^- \left( \frac12 \phi_m + \DL_{m,mm} \phi_m - \frac{\sigma_0}{\sigma_m} \SL_{m,mm} \psi_m \right)
		= \eta_m^- \traceDme \uinc, \\
		\nu_m^+ \left( \frac12 \psi_m + \sum_{n=1}^\ell \left(\HS_{0,mn} \phi_n + \AD_{0,mn} \psi_n\right) \right) \\
		\quad + \nu_m^- \left( \frac12 \psi_m - \frac{\sigma_m}{\sigma_0} \HS_{m,mm} \phi_m - \AD_{m,mm} \psi_m \right)
		= \nu_m^+ \traceNme \uinc
	\end{cases}
\end{equation}
for $m=1,2,\dots,\ell$.
The \emph{combined domain formulation} is given by
\begin{equation}
	\label{eq:formulation:combined:domain}
	\begin{cases}
		\eta_m^+ \left( \frac12 \phi_m + \sum_{n=1}^\ell \left(-\DL_{0,mn} \phi_n + \SL_{0,mn} \psi_n\right) \right) \\
		\quad + \nu_m^+ \left( \frac12 \psi_m + \sum_{n=1}^\ell \left(\HS_{0,mn} \phi_n + \AD_{0,mn} \psi_n\right) \right) \\
		\quad = \eta_m^- \traceDme \uinc + \nu_m^+ \traceNme \uinc, \\
		\eta_m^- \left( \frac12 \phi_m + \DL_{m,mm} \phi_m - \frac{\sigma_0}{\sigma_m} \SL_{m,mm} \psi_m \right) \\
		\quad + \nu_m^- \left( \frac12 \psi_m - \frac{\sigma_m}{\sigma_0} \HS_{m,mm} \phi_m - \AD_{m,mm} \psi_m \right)
		= 0
	\end{cases}
\end{equation}
for $m=1,2,\dots,\ell$.
The \emph{combined mixed formulation} is given by
\begin{equation}
	\label{eq:formulation:combined:mixed}
	\begin{cases}
		\eta_m^+ \left( \frac12 \phi_m + \sum_{n=1}^\ell \left(-\DL_{0,mn} \phi_n + \SL_{0,mn} \psi_n\right) \right) \\
		\quad + \nu_m^- \left( \frac12 \psi_m - \frac{\sigma_m}{\sigma_0} \HS_{m,mm} \phi_m - \AD_{m,mm} \psi_m \right)
		= \eta_m^- \traceDme \uinc, \\
		\nu_m^+ \left( \frac12 \psi_m + \sum_{n=1}^\ell \left(\HS_{0,mn} \phi_n + \AD_{0,mn} \psi_n\right) \right) \\
		\quad + \eta_m^- \left( \frac12 \phi_m + \DL_{m,mm} \phi_m - \frac{\sigma_0}{\sigma_m} \SL_{m,mm} \psi_m \right)
		= \nu_m^+ \traceNme \uinc
	\end{cases}
\end{equation}
for $m=1,2,\dots,\ell$.

\paragraph{Remarks}
The first references to single-trace formulations were for single objects only and were quickly extended to multiple reflection~\cite{arvas1986field}. In acoustics, the Müller formulation has also been called the Burton-Miller formulation for penetrable domains~\cite{wu2012fast}. In electromagnetics, the combined single-trace formulations are also known as \emph{combined field integral equations} (CFIE)~\cite{mautz1977electromagnetic, harrington1989boundary, rao1990field}. Finally, these formulations have also been used to solve diffusion equations~\cite{sikora2006diffuse} and Poisson-Boltzmann systems~\cite{juffer1991electric, bardhan2009numerical}.

\subsubsection{Multiple-traces formulations}

The multiple-traces formulations use the four surface potentials
\begin{equation}
	\phi_m^+ = \traceDme\utot,
	\ \psi_m^+ = \traceNme\utot,
	\ \phi_m^- = \traceDmi\utot,
	\ \text{and}
	\ \psi_m^- = \traceNmi\utot.
\end{equation}
Then, the interface conditions are used to convert the potentials related to the identity operators that arise in the traces of the direct representation formulas~\eqref{eq:representation:exterior:direct}--\eqref{eq:representation:interior:direct}. This results in
\begin{equation}
	\begin{cases}
		\frac12 \begin{bmatrix} \ID_m & 0 \\ 0 & \frac{\sigma_m}{\sigma_0} \ID_m \end{bmatrix} \begin{bmatrix} \phi_m^- \\ \psi_m^- \end{bmatrix}
		+ \sum_{n=1}^\ell \CA_{0,mn} \begin{bmatrix} \phi_n^+ \\ \psi_n^+ \end{bmatrix}
		= \begin{bmatrix} \traceDme \uinc \\ \traceNme \uinc \end{bmatrix}, \\
		-\frac12 \begin{bmatrix} \ID_m & 0 \\ 0 & \frac{\sigma_0}{\sigma_m} \ID_m \end{bmatrix} \begin{bmatrix} \phi_m^+ \\ \psi_m^+ \end{bmatrix}
		+ \CA_m \begin{bmatrix} \phi_m^- \\ \psi_m^- \end{bmatrix}
		= \begin{bmatrix} 0 \\ 0 \end{bmatrix}
	\end{cases}
	\label{eq:formulation:mtf}
\end{equation}
for $m=1,2,\dots,\ell$ which is called the \emph{multiple-traces formulation}~\cite{hiptmair2012multiple, peng2012computations}. Other versions of multiple-traces formulations include global interconnections~\cite{claeys2015integral} and combined fields~\cite{claeys2019introduction}.

\subsubsection{Auxiliary field formulations}

Previously, the interior fields took the value of the pressure field in one subdomain and zero outside. Here, let us consider interior fields defined as
\begin{align}
	p_m &= \begin{cases} \utot & \text{in } \Omega_m, \\ \hat{p}_m & \text{in } \Omega_n \text{ for } n=0,1,2,\dots,\ell \text{ and } n \ne m; \end{cases} \label{formulation:field:interior:auxiliary}
\end{align}
for $m=1,2,\dots,\ell$ where $\hat{p}_m$ is an unknown auxiliary field exterior to subdomain~$m$. Now, the direct representation formula~\eqref{eq:representation:interior:direct} reads
\begin{align}
	p_m &= \SLP_m \hat\psi_m - \DLP_m \hat\phi_m
	\label{eq:representation:interior:auxiliary}
\end{align}
for $m = 1,2,\dots,\ell$ where the auxiliary potentials have to satisfy
\begin{align}
	\hat\phi_m &= \traceDmi \utot - \traceDme \hat{p}_m = \traceDme \left(\utot - \hat{p}_m\right), \label{eq:auxiliary:phi} \\
	\hat\psi_m &= \traceNmi \utot - \traceNme \hat{p}_m = \traceNme \left(\frac{\sigma_0}{\sigma_m} \utot - \hat{p}_m\right) \label{eq:auxiliary:psi}
\end{align}
for $m=1,2,\dots,\ell$ because of the jump relations for boundary integral operators~\cite{steinbach2008numerical}.
The exterior traces of the auxiliary representation formula~\eqref{eq:representation:interior:auxiliary} yield the auxiliary Calderón system
\begin{align}
	\begin{bmatrix} \traceDme \hat{p}_m \\ \traceNme \hat{p}_m \end{bmatrix}
	&= \left(-\frac12\ID_m + A_m \right)
	\begin{bmatrix} \hat\phi_m \\ \hat\psi_m \end{bmatrix}
\end{align}
for $m=1,2,\dots,\ell$.
Since the auxiliary fields are arbitrary, one can impose either $\traceDme \hat{p}_m = \traceDme \utot$ or $\traceNme \hat{p}_m = \frac{\sigma_0}{\sigma_m} \traceNme \utot$, but not both. Then, the auxiliary field cannot be zero, which prevents obtaining the standard single-trace formulations. Furthermore, these choices lead to either
\begin{equation}
	\label{eq:auxiliarytraces:psi}
	\begin{cases}
		\hat\phi_m = 0, \\
		\traceDme \utot = \SL_m \hat\psi_m, \\
		\traceNme \utot = \frac{\sigma_m}{\sigma_0} \left( \frac12\ID_m + \AD_m \right) \hat\psi_m
	\end{cases}
\end{equation}
for $m=1,2,\dots,\ell$; or
\begin{equation}
	\label{eq:auxiliarytraces:phi}
	\begin{cases}
		\hat\psi_m = 0, \\
		\traceDme \utot = \left( \frac12\ID_m - \DL_m \right) \hat\phi_m, \\
		\traceNme \utot = \frac{\sigma_m}{\sigma_0} \HS_m \hat\phi_m
	\end{cases}
\end{equation}
for $m=1,2,\dots,\ell$, respectively.

Now, boundary integral formulations can be designed by taking the exterior traces of the direct exterior representation formula~\eqref{eq:representation:exterior:direct} and substitute the auxiliary traces~\eqref{eq:auxiliarytraces:psi} to obtain
\begin{align}
	&\frac12 \SL_m \hat\psi_m + \sum_{n=1}^\ell \left( -\DL_{0,mn} \SL_n + \frac{\sigma_n}{\sigma_0} \SL_{0,mn} \left( \frac12\ID_n + \AD_n \right) \right) \hat\psi_n
	= \traceDm \uinc, \label{eq:formulation:auxiliary:dir:psi} \\
	&\frac12 \frac{\sigma_m}{\sigma_0} \left( \frac12\ID_m + \AD_m \right) \hat\psi_m + \sum_{n=1}^\ell \left( \HS_{0,mn} \SL_n + \frac{\sigma_n}{\sigma_0} \AD_{0,mn} \left( \frac12\ID_n + \AD_n \right) \right) \hat\psi_n \nonumber \\
	&\quad = \traceNm \uinc \label{eq:formulation:auxiliary:neu:psi}
\end{align}
for $m=1,2,\dots,\ell$ or, alternatively, substitute the auxiliary traces~\eqref{eq:auxiliarytraces:phi} to obtain
\begin{align}
	&\frac12 \left( \frac12\ID_m - \DL_m \right) \hat\phi_m + \sum_{n=1}^\ell \left( -\DL_{0,mn} \left( \frac12\ID_n - \DL_n \right) + \frac{\sigma_n}{\sigma_0} \SL_{0,mn} \HS_n \right) \hat\phi_n \nonumber \\
	&\quad = \traceDme \uinc, \label{eq:formulation:auxiliary:dir:phi} \\
	&\frac12 \frac{\sigma_m}{\sigma_0} \HS_m \hat\phi_m + \sum_{n=1}^\ell \left( \HS_{0,mn} \left( \frac12\ID_n - \DL_n \right) + \frac{\sigma_n}{\sigma_0} \AD_{0,mn} \HS_n \right) \hat\phi_n \nonumber \\
	&\quad = \traceNme \uinc \label{eq:formulation:auxiliary:neu:phi}
\end{align}
for $m=1,2,\dots,\ell$.
The four boundary integral formulations~\eqref{eq:formulation:auxiliary:dir:psi}--\eqref{eq:formulation:auxiliary:neu:phi} are called \emph{auxiliary field formulations} and have only one unknown surface potential at each interface. Notice that linear combinations of formulations~\eqref{eq:formulation:auxiliary:dir:psi} and~\eqref{eq:formulation:auxiliary:neu:psi} or formulations~\eqref{eq:formulation:auxiliary:dir:phi} and~\eqref{eq:formulation:auxiliary:neu:phi} could be used as well.

These formulations were introduced for time-dependent problems~\cite{marx1982single, marx1984integral}, are also known as single-source formulations in electromagnetics~\cite{glisson1984integral, mautz1989stable}, and can also be derived using an indirect approach with a specific set of potentials~\cite{kleinman1988single}.

\subsection{Indirect boundary integral formulations}
\label{sec:formulation:indirect}

Any solution of the Helmholtz system can be represented by an indirect representation of the field in terms of a single surface potential~\cite{steinbach2008numerical}, which is not necessarily the trace of the pressure field. This leads to the \emph{single-potential formulations} (SPF) while mixing direct and indirect representations for the exterior and interior result in the \emph{mixed-potential formulations} (MPF).

\subsubsection{Single potential formulations}
\label{sec:indirect:singlepotential}

The fields are represented by a single surface potential at each interface as
\begin{align}
	\usca &= -\sum_{n=1}^\ell \SLP_{0,n} \psi_m^+ & \text{or}\quad
	\usca &= \sum_{n=1}^\ell \DLP_{0,n} \phi_m^+ & \text{in } \Omega_0,
	\label{eq:representation:exterior:single} \\
	\utot &= \SLP_m \psi_m^- & \text{or}\quad
	\utot &= -\DLP_m \phi_m^- & \text{in } \Omega_m
	\label{eq:representation:interior:single}
\end{align}
for $m=1,2,\dots,\ell$, and surface potentials $\phi_m^\pm$ and $\psi_m^\pm$ which are not necessarily the traces of the pressure field.
Substituting the single-layer representation into the interface conditions~\eqref{eq:helmholtz} yields
\begin{align}
	\begin{cases}
		\SL_m \psi_m^- - \sum_{n=1}^\ell \SL_{0,mn} \psi_n^+ = \traceDme \uinc, \\
		\frac{\sigma_m}{\sigma_0} \left(\frac12\ID_m + \AD_m\right) \psi_m^- + \frac12\ID_m \psi_m^+ - \sum_{n=1}^\ell \AD_{0,mn} \psi_n^+ = \traceNme \uinc
	\end{cases}
	\label{eq:formulation:singlepotential:slpext:slpint}
\end{align}
for $m=1,2,\dots,\ell$.
Alternatively, taking the double-layer representation yields
\begin{align}
	\begin{cases}
		\left(\frac12\ID_m - \DL_m\right) \phi_m^- + \frac12\ID_m \phi_m^+ + \sum_{n=1}^\ell \DL_{0,mn} \phi_n^+ = \traceDme \uinc, \\
		\frac{\sigma_m}{\sigma_0} \HS_m \phi_m^- - \sum_{n=1}^\ell \HS_{0,mn} \phi_n^+ = \traceNme \uinc
	\end{cases}
	\label{eq:formulation:singlepotential:dlpext:dlpint}
\end{align}
for $m=1,2,\dots,\ell$.
Two other boundary integral formulations can be designed by mixing the indirect single-layer and double-layer representations, that is,
\begin{align}
	\begin{cases}
		\left(\frac12\ID_m - \DL_m\right) \phi_m^- - \sum_{n=1}^\ell \SL_{0,mn} \psi_n^+ = \traceDme \uinc, \\
		\frac{\sigma_m}{\sigma_0} \HS_m \phi_m^- + \frac12\ID_m \psi_m^+ - \sum_{n=1}^\ell \AD_{0,mn} \psi_n^+ = \traceNme \uinc
	\end{cases}
	\label{eq:formulation:singlepotential:slpext:dlpint}
\end{align}
for $m=1,2,\dots,\ell$ and
\begin{align}
	\begin{cases}
		\SL_m \psi_m^- + \frac12\ID_m \phi_m^+ + \sum_{n=1}^\ell \DL_{0,mn} \phi_n^+ = \traceDme \uinc, \\
		\frac{\sigma_m}{\sigma_0} \left(\frac12\ID_m + \AD_m\right) \psi_m^- - \sum_{n=1}^\ell \HS_{0,mn} \phi_n^+ = \traceNme \uinc
	\end{cases}
	\label{eq:formulation:singlepotential:dlpext:slpint}
\end{align}
for $m=1,2,\dots,\ell$.
Furthermore, one could also consider different interior potentials at different interfaces.
Analogous versions for Maxwell's equations have only an electric or magnetic surface currents as unknown potentials and are called the \emph{electric} and \emph{magnetic current formulations}~\cite{harrington1989boundary, yla2011calderon}.

\subsubsection{Mixed potential formulations}

Let us consider a mix of single-potential representation for the exterior and a direct representation for the interior fields, that is,
\begin{align}
	\usca &= -\sum_{n=1}^\ell \SLP_{0,n} \psi_m^+ \quad\text{or}\quad
	\usca = \sum_{n=1}^\ell \DLP_{0,n} \phi_m^+ & \text{in } \Omega_0,
	\label{eq:representation:exterior:mpf} \\
	\utot &= \SLP_m (\traceNmi \utot) - \DLP_m (\traceDmi \utot) & \text{in } \Omega_m
	\label{eq:representation:interior:mpf}
\end{align}
for $m=1,2,\dots,\ell$, and surface potentials $\phi_m^+$ and $\psi_m^+$ which are not necessarily the traces of the pressure field.
Since a direct formulation is used for the interior field, one can use the interior NtD and DtN maps~\eqref{eq:ntd}--\eqref{eq:dtn}. Then, the interface conditions~\eqref{eq:helmholtz} yield
\begin{align}
	\traceDme \utot &= \traceDmi \utot = \NtDmi \traceNmi \utot = \frac{\sigma_0}{\sigma_m} \NtDmi \traceNme \utot, \\
	\traceNme \utot &= \frac{\sigma_m}{\sigma_0} \traceNmi \utot = \frac{\sigma_m}{\sigma_0} \DtNmi \traceDmi \utot = \frac{\sigma_m}{\sigma_0} \DtNmi \traceDme \utot.
\end{align}
Now, boundary integral formulations can be designed by taking traces of the exterior field and eliminating one of the traces via the NtD or DtN maps.
Specifically, taking the exterior single-layer representation and eliminating the exterior Neumann trace results in
\begin{align}
	\begin{cases}
		\traceDme \utot - \sum_{n=1}^\ell \SL_{0,mn} \psi_{0,n} = \traceDme \uinc, \\
		\frac{\sigma_m}{\sigma_0} \DtNmi (\traceDme \utot) + \frac12 \psi_{0,m} - \sum_{n=1}^\ell \AD_{0,mn} \psi_{0,n} = \traceNme \uinc
	\end{cases}
	\label{eq:formulation:mixedpotential:slp:dtn}
\end{align}
for $m=1,2,\dots,\ell$ and the unknowns $\psi_{0,m}$ and $\traceDme \utot$,
while eliminating the Dirichlet trace results in
\begin{align}
	\begin{cases}
		\frac{\sigma_0}{\sigma_m} \NtDmi (\traceNme \utot) - \sum_{n=1}^\ell \SL_{0,mn} \psi_{0,n} = \traceDme \uinc, \\
		\traceNme \utot + \frac12 \psi_{0,m} - \sum_{n=1}^\ell \AD_{0,mn} \psi_{0,n} = \traceNme \uinc
	\end{cases}
	\label{eq:formulation:mixedpotential:slp:ntd}
\end{align}
for $m=1,2,\dots,\ell$ and the unknowns $\psi_{0,m}$ and $\traceNme \utot$.
Alternatively, taking the exterior double-layer representation and eliminating the exterior Neumann trace results in
\begin{equation}
	\begin{cases}
		\traceDme \utot + \frac12 \phi_{0,m} + \sum_{n=1}^\ell \DL_{0,mn} \phi_{0,n} = \traceDme \uinc \\
		\frac{\sigma_m}{\sigma_0} \DtNmi (\traceDme \utot) - \sum_{n=1}^\ell \HS_{0,mn} \phi_{0,n} = \traceNme \uinc.
	\end{cases}
	\label{eq:formulation:mixedpotential:dlp:dtn}
\end{equation}
for $m=1,2,\dots,\ell$ and the unknowns $\phi_{0,m}$ and $\traceDme \utot$,
while eliminating the Dirichlet trace results in
\begin{align}
	\begin{cases}
		\frac{\sigma_0}{\sigma_m} \NtDmi (\traceNme \utot) + \frac12 \phi_{0,m} + \sum_{n=1}^\ell \DL_{0,mn} \phi_{0,n} = \traceDme \uinc \\
		\traceNme \utot - \sum_{n=1}^\ell \HS_{0,mn} \phi_{0,n} = \traceNme \uinc.
	\end{cases}
	\label{eq:formulation:mixedpotential:dlp:ntd}
\end{align}
for $m=1,2,\dots,\ell$ and the unknowns $\phi_{0,m}$ and $\traceNme \utot$.

These formulations include NtD and DtN maps that have no closed-form expressions for general surfaces. Hence, the equations need to be multiplied from the left by the correct boundary integral operators according to the definitions of the NtD and DtN maps~\eqref{eq:ntd}--\eqref{eq:dtn}. For example, the formulation~\eqref{eq:formulation:mixedpotential:slp:dtn} can be written as
\begin{align}
	\begin{cases}
		\traceDme \utot - \sum_{n=1}^\ell \SL_{0,mn} \psi_{0,n} = \traceDme \uinc, \\
		\frac{\sigma_m}{\sigma_0} \HS_m (\traceDme \utot) + \left(\frac12\ID_m - \AD_m\right) \left( \frac12 \psi_{0,m} - \sum_{n=1}^\ell \AD_{0,mn} \psi_{0,n} \right) \\
		\quad = \left(\frac12\ID_m - \AD_m\right) \traceNme \uinc
	\end{cases}
	\label{eq:formulation:mixedpotential:slp:dtn:ad}
\end{align}
or
\begin{align}
	\begin{cases}
		\traceDme \utot - \sum_{n=1}^\ell \SL_{0,mn} \psi_{0,n} = \traceDme \uinc, \\
		\frac{\sigma_m}{\sigma_0} \left(\frac12\ID_m + \DL_m\right) (\traceDme \utot) + \SL_m \left( \frac12 \psi_{0,m} - \sum_{n=1}^\ell \AD_{0,mn} \psi_{0,n} \right) \\
		\quad = \SL_m \traceNme \uinc
	\end{cases}
	\label{eq:formulation:mixedpotential:slp:dtn:sl}
\end{align}
for $m=1,2,\dots,\ell$. This procedure yields eight different \emph{mixed potential formulations}.
Notice that reversing the direct and indirect formulation is complicated in the case of multiple reflection since exterior NtD and DtN maps will not be local anymore.
Formulations based on the same design principles are called \emph{single-source formulations} in electromagnetics~\cite{gossye2018calderon, gossye2019electromagnetic} and are well conditioned for high-contrast media~\cite{wout2022highcontrast}.

\subsection{Other boundary integral formulations}

The list of boundary integral formulations presented above is not exhaustive.
Firstly, indirect formulations can be designed with explicit relations between surface potentials, such as the combined-source formulations in electromagnetics~\cite{harrington1989boundary}, the Brakhage-Werner formulation~\cite{rapun2006indirect} and regularised formulations~\cite{boubendir2015integral}. Secondly, domain decomposition techniques result in global multiple-traces formulations~\cite{claeys2013multi}, symmetric mortar element formulations for five unknown potentials at each interface~\cite{laliena2009symmetric}, or specific coupling conditions for independent subdomain formulations~\cite{langer2003boundary, wu2008multi, peng2012nonconformal}. Thirdly, electromagnetic formulations can be designed with respect to the scalar and vector potentials, in addition to the fields (augmented formulations~\cite{xia2016enhanced}) or as replacements of the fields ($A$-$\phi$ formulations~\cite{chew2014vector}). We do not claim completeness of the formulations mentioned in this study since there is an abundance of literature on the topic and the design of novel formulations is still actively pursued.

\subsection{Numerical discretisation}
\label{sec:discretisation}

The boundary integral operators can be discretised with numerical algorithms such as collocation and Galerkin methods. Here, a Galerkin method is used with continuous piecewise linear (P1) functions as test and basis functions in the weak formulation~\cite{smigaj2015solving}. Given a triangulation of the material interface, the P1 functions have a value of one on a specific node, are zero on all other nodes in the mesh, and have a continuous linear approximation on each triangular element.

\section{Preconditioning}
\label{sec:preconditioning}

The discretised boundary integral formulations are a dense system of linear equations that need to be solved with either direct factorisation~\cite{bremer2015high} or iterative Krylov methods~\cite{marburg2003performance}. For large-scale simulations, the GMRES algorithm~\cite{saad1986gmres} is often the preferred technique, where the dense matrix arithmetic is accelerated with the fast multipole method~\cite{greengard1987fast} or hierarchical matrix compression~\cite{hackbusch2015hierarchical}. Furthermore, preconditioning of the linear system is often essential to limit the number of GMRES iterations to reach a predefined accuracy. Since algebraic preconditioners such as ILU require explicit access to the matrix~\cite{schneider2003performance}, they are cumbersome to implement in conjunction with acceleration techniques~\cite{sakuma2008fast}. Differently, operator preconditioning is based on boundary integral operators that are discretised separately to the model formulation and can, therefore, readily be combined with accelerators~\cite{darbas2013combining}. Moreover, the effectiveness of operator preconditioners is justified by functional analysis of the boundary integral operators~\cite{antoine2008integral}.

Given a linear operator $Q:V \to W$ that maps from function space $V$ to $W$, a precondioner with the mapping property $R:W \to V$ yields a system $RQ:V \to V$ that is typically well conditioned~\cite{hiptmair2006operator, kirby2010functional}. This observation is the basis of so-called operator preconditioning. In the case of weak formulations, the discretised operators satisfy $Q_h:V_h \to W'_h$ and $R_h:W_h \to V'_h$ where the subscript $h$ denotes a finite-dimensional subspace and the prime denotes the dual space. Then, additional mass matrices achieve the desired mapping property of
\begin{equation*}
	M_2^{-1} R_h M_1^{-1} Q_h:V_h \to V_h
\end{equation*}
for $M_1:W_h \to W_h'$ and $M_2:V_h \to V_h'$ discretised identity operators~\cite{betcke2020product}. In this study, mass-matrix preconditioning will always be used. Moreover, the operator products for the preconditioned formulations are not explicitly calculated but separate matrix-vector multiplications are performed at each iteration of the iterative linear solver.

Two of the most effective preconditioning strategies are Calderón and OSRC preconditioning. Table~\ref{table:formulations} summarises the feasible combinations of preconditioner and model formulation, along with characteristics of the preconditioned boundary integral formulations.

\begin{table}[!ht]
	\caption{Preconditioned boundary integral formulations for acoustic transmission throught multiple objects. The abbreviations are Calderón projection (CP), on-surface radiation conditions (OSRC), opposite-order (OO), unknown surface potentials (sp), dense boundary integral operators (BIO), dense matrix-vector multiplications (mv), and $\ell$ is the number of bounded domains.}
	\label{table:formulations}
	\bigskip
	\centering
	\begin{tabular}{lcccccc}
		formulation & CP & OSRC/OO & \#sp & \#BIO & \#mv \\
		\hline
		Dirichlet \eqref{eq:formulation:dirichlet}
		& - & $\checkmark$ & $2\ell$ & $2\ell+2\ell^2$ & $2\ell+2\ell^2$ \\
		Neumann \eqref{eq:formulation:neumann}
		& - & $\checkmark$ & $2\ell$ & $2\ell+2\ell^2$ & $2\ell+2\ell^2$ \\
		PMCHWT \eqref{eq:formulation:pmchwt}
		& $\checkmark$ & $\checkmark$ & $2\ell$ & $4\ell+4\ell^2$ & $4\ell+4\ell^2$ \\
		Müller \eqref{eq:formulation:muller}
		& - & - & $2\ell$ & $4\ell+4\ell^2$ & $4\ell+4\ell^2$ \\
		Combined field \eqref{eq:formulation:combined:trace}--\eqref{eq:formulation:combined:mixed}
		& - & $\checkmark$ & $2\ell$ & $4\ell+4\ell^2$ & $4\ell+4\ell^2$ \\
		Multiple traces \eqref{eq:formulation:mtf}
		& $\checkmark$ & $\checkmark$ & $4\ell$ & $4\ell+4\ell^2$ & $4\ell+4\ell^2$ \\
		Auxiliary field \eqref{eq:formulation:auxiliary:dir:psi}--\eqref{eq:formulation:auxiliary:neu:phi}
		& - & $\checkmark$ & $\ell$ & $2\ell+2\ell^2$ & $\ell+3\ell^2$ \\
		Single potential \eqref{eq:formulation:singlepotential:slpext:slpint}--\eqref{eq:formulation:singlepotential:dlpext:slpint}
		& - & $\checkmark$ & $2\ell$ & $2\ell+2\ell^2$ & $2\ell+2\ell^2$ \\
		Mixed potential \eqref{eq:formulation:mixedpotential:slp:dtn}--\eqref{eq:formulation:mixedpotential:dlp:ntd}
		& - & $\checkmark$ & $2\ell$ & $2\ell+2\ell^2$ & $\ell+3\ell^2$ \\
	\end{tabular}
\end{table}

\subsection{Calderón preconditioning}

The family of \emph{Calderón preconditioners} are designed with information from the Calderón identities introduced in Section~\ref{sec:calderonidentities}. These preconditioners are linear operators with dense blocks and are mainly effective at moderate frequency ranges.

\subsubsection{Projection-based preconditioning}

The Calderón operator is a projection, specifically, $\CA_m^2 = \frac14 \bar\ID_m$ and $\widehat{\CA}_m^2 = \frac14 \bar\ID_m$. Hence, the Calderón operator is a perfect preconditioner for itself. This property is the design principle behind Calderón preconditioning of the PMCHWT and MTF formulations.

The PMCHWT formulation~\eqref{eq:formulation:pmchwt} involves sums of interior and exterior Calderón operators. Then, Calderón preconditioning is justified with the following observations:
\begin{align}
	&\left(A_0 + \widehat{A}_1\right)^2
	= \frac14 \bar\ID_0 + \frac14 \bar\ID_1 + A_0 \widehat{A}_1 + \widehat{A}_1 A_0, \\
	&A_0 \left(A_0 + \widehat{A}_1\right)
	= \frac14 \ID_0 + A_0 \widehat{A}_1, \\
	&\widehat{A}_1 \left(A_0 + \widehat{A}_1\right)
	= \frac14 \ID_1 + \widehat{A}_1 A_0.
\end{align}
These are well-conditioned operators since the products of Calderón operators are compact~\cite{antoine2008integral, claeys2013multi, boubendir2015integral}. The full version incurs more computation time than the exterior and interior preconditioners but tends to be better conditioned since it has a single spectral accumulation point~\cite{wout2022highcontrast}. None of the preconditioners require additional storage since the operators are already present in the model anyway. This advantage is lost when a different wavenumber or numerical parameters are chosen for Calderón preconditioning~\cite{yan2010comparative, boubendir2015integral}.

The boundary integral formulations for transmission through multiple objects result in a linear system with blocks associated to each material interface. Calderón preconditioning based on the full system is effective but requires an overhead of $\ell + \ell^2$ matrix-vector multiplications of individual Calderón matrices. A more efficient alternative is a diagonal block preconditioner based on a single Calderón operator at each interface. This approach does not incorporate multiple reflection in the preconditioner but has a superior computational complexity of only $\ell$ matrix-vector multiplications in the preconditioner step of GMRES.

\subsubsection{Opposite-order preconditioning}
\label{sec:preconditioning:oppositeorder}

A corollary of operator preconditioning for Sobolev spaces is that the preconditioner needs to be of opposite order compared to the model. Hence, single-layer and hypersingular boundary integral operators are good candidates for preconditioning respectively the hypersingular and single boundary integral operator~\cite{steinbach1998construction}. The efficiency of this \emph{opposite-order preconditioning} is also justified by the Calderón identities~\eqref{eq:slhs}--\eqref{eq:hssl}. For instance, the single-potential formulation~\eqref{eq:formulation:singlepotential:dlpext:dlpint} can be preconditioned as
\begin{align}
	\begin{bmatrix}
		\ID_1 & 0 \\
		0 & -\SL_0
	\end{bmatrix}
	\begin{bmatrix}
		\frac12\ID_1 - \DL_1 & \frac12\ID_0 + \DL_0 \\
		\frac{\sigma_1}{\sigma_0} \HS_1 & -\HS_0 
	\end{bmatrix}
	\begin{bmatrix} \phi_1 \\ \phi_0 \end{bmatrix}
	= \begin{bmatrix}
		\ID_1 & 0 \\
		0 & -\SL_0
	\end{bmatrix}
	\begin{bmatrix} \traceDme \uinc \\ \traceNme \uinc \end{bmatrix}
\end{align}
for a single domain and with direct extensions to multiple domains.

\subsection{OSRC preconditioning}

The DtN and NtD maps are good candidates as opposite-order preconditioners for the hypersingular and single-layer operators since this combination yields a second-kind operator, as can be seen in Eqns.~\eqref{eq:ntd:prec}--\eqref{eq:dtn:prec}. However, directly using these maps is not feasible since no closed-form expressions are available for general surfaces. Hence, approximations need to be used, of which the on-surface radiation conditions are among the most efficient ones~\cite{moore1988theory}. The OSRC preconditioners are local operators and are especially accurate at high frequencies~\cite{antoine2008advances}.
They are defined by the pseudo-differential operators
\begin{align}
	\OsrcNtDmi &= \frac1{\imath k_m} \left(\ID_m + \frac{\Delta_{\Gamma_m}}{k_{m,\epsilon}^2}\right)^{-\frac12}, \\
	\OsrcDtNmi &= \imath k_m \left(\ID_m + \frac{\Delta_{\Gamma_m}}{k_{m,\epsilon}^2}\right)^{\frac12}
\end{align}
where $\Delta_{\Gamma_m}$ denotes the Laplace-Beltrami operator on surface~$\Gamma_m$ and $k_{m,\epsilon} = k_m (1+\imath\epsilon)$ a damped wavenumber for $\epsilon>0$. Setting the hyperparameters is often based on optimal choices for a single spherical geometry~\cite{antoine2005alternative}. Specifically, $\epsilon = 0.4 (k_m R_m)^{-\frac23}$ where $R_m$ denotes the radius of the object~$\Omega_m$. The square-root operation will be approximated with a truncated Padé series expansion that reduces the operator into a set of surface Helmholtz equations with complex-valued wavenumbers~\cite{darbas2013combining}. Since these are local boundary integral operators, the resulting matrices are sparse, the inversion of which is performed by calculating the sparse LU-factorisation once, and stored for use in each iteration of the linear solver.

As an example of OSRC preconditioning, the PMCHWT formulation for a single object becomes
\begin{align*}
	\begin{bmatrix}
		0 & V_{\mathrm{NtD},1}^- \\
		V_{\mathrm{DtN},1}^- & 0
	\end{bmatrix}
	\begin{bmatrix}
		-\DL_0 - \DL_1 & \SL_0 + \frac{\sigma_0}{\sigma_1} \SL_1 \\
		\HS_0 + \frac{\sigma_1}{\sigma_0} \HS_1 & \AD_0 + \AD_1
	\end{bmatrix}
	\begin{bmatrix} \phi \\ \psi \end{bmatrix}
	= \begin{bmatrix}
		0 & V_{\mathrm{NtD},1}^- \\
		V_{\mathrm{DtN},1}^- & 0
	\end{bmatrix}
	\begin{bmatrix} \traceDe \uinc \\ \traceNe \uinc \end{bmatrix}
\end{align*}
which is equivalent to a block-diagonal preconditioner for a permuted PMCHWT formulation~\cite{haqshenas2020fast} and with direct extensions to multiple domains~\cite{wout2022pmchwt}. Furthermore, the OSRC operators can be used as a combination parameter for the combined single-trace formulations~\eqref{eq:formulation:combined:trace}--\eqref{eq:formulation:combined:mixed} as well~\cite{darbas2013combining, wout2015fast, betcke2017computationally}.

\section{Benchmarking}
\label{sec:results}

The previous sections presented five families of boundary integral formulations and three preconditioning strategies. This section assesses the computational characteristics of these preconditioned formulations through an extensive benchmarking exercise. The objective is to compare different preconditioned formulation, rather than different mathematical models~\cite{wout2021proximity}, physical experiments, or software packages~\cite{hornikx2015platform}.

\subsection{Parameter selection}

Any benchmarking excercise needs to limit the parameter space to a feasible size to perform the computational simulations. Here, the choices to set the parameters will be explained.

\paragraph{Wave propagation}
The incident wave field is a plane wave with driving frequency~$f$ in Hz. The physical parameters are chosen from materials commonly found in biomedical engineering, with a linear frequency power law model for attenuation~\cite{hamilton1998nonlinear}. That is, $\sigma_m = 1/\rho_m$ for $\rho$ the mass density in kg~m$^{-3}$ and $k_m= 2\pi f/c_m + \iota \alpha_m (f \cdot 10^{-6})^{b_m}$ where $c$ denotes the wavespeed in m~s$^{-1}$, $\alpha$ the attenuation coefficient in Np~m$^{-1}$~Hz$^{-1}$ and $b$ an exponent. See Table~\ref{table:parameters:physical} for the values.

\begin{table}[!ht]
	\caption{Physical parameters of wave propagation through different materials~\cite{duck1990physical,itis2018}.}
	\label{table:parameters:physical}
	\centering
	\begin{tabular}{lrrrr}
		\hline\hline
		material & $\rho$ & $c$ & $\alpha$ & $b$ \\
		\hline
		water & 1000 & 1500 & 0.015 & 2 \\
		fat & 917 & 1412 & 9.334 & 1 \\
		bone & 1912 & 4080 & 47.20 & 1 \\
		\hline\hline
	\end{tabular}
\end{table}

\paragraph{Numerical discretisation}
As explained in Section~\ref{sec:discretisation}, a Galerkin method with P1 elements is used as numerical discretisation. The benchmarks will be performed with dense matrix algebra, even though all preconditioned formulations are amenable to acceleration schemes. The main reason is to limit the parameter space: fast multipole and matrix compression algorithms often require expert choices for numerical parameters to obtain fast implementations. Standard quadrature rules are used with three points per triangle, which are increased to five for near interactions, and a singularity-aware scheme for the self interactions.

\paragraph{Meshes}
The triangular surface meshes are generated with the open-source library Gmsh~\cite{geuzaine2009gmsh}. The triangular elements are flat, no curved elements are used~\cite{zhang2020dual}. The edges of the triangles have a length of at most $h = \lambda_\mathrm{min}/n_h$ where $\lambda_\mathrm{min}$ is the minimum of the wavelengths in the interior and exterior region to each interface, and $n_h$ is a fixed number. Here, at least four elements per wavelength are used ($n_h = 4$). Even though this is on the lower side of common choices for the mesh width\cite{marburg2002six}, a Galerkin method with P1 elements is sufficiently accurate at a sphere~\cite{haqshenas2020fast}. For comparison, the benchmarks at 1~MHz were performed with eight elements per wavelength as well, without yielding different conclusions for this study.

\paragraph{Geometry}
The objects are spheres with a radius of 5~mm and acoustic properties of either fat or bone, with water being the exterior medium. In the case of two spheres, one is made of fat and the other is bone, with a distance of 35~mm between the two centers. Table~\ref{table:parameters:numerical} summarizes the geometrical details.

\begin{table}[!ht]
	\caption{Numerical parameters of the geometries, with $D$ the domain size and $\lambda$ the minimum wavelength of the materials.}
	\label{table:parameters:numerical}
	\centering
	\begin{tabular}{rrrrrrr}
		\hline\hline
		& \multicolumn{2}{c}{water-fat} & \multicolumn{2}{c}{water-bone} & \multicolumn{2}{c}{water-fat-bone} \\
		frequency & $D/\lambda$ & \#nodes & $D/\lambda$ & \#nodes & $D/\lambda$ & \#nodes \\
		\hline
		1 MHz   &  7.08 &  3246 &  6.67 &  2777 & 31.87 & 6498 \\
		1.5 MHz & 10.62 &  7086 & 10.00 &  6302 \\
		2 MHz   & 14.16 & 12377 & 13.33 & 10782 \\
		\hline\hline
	\end{tabular}
\end{table}

\paragraph{Preconditioned linear solver}
All linear systems were solved with GMRES~\cite{saad1986gmres}, with a termination criterion of $\epsilon = 10^{-5}$, a maximum of 1000 iterations, and without restart.
The Calderón preconditioners reuse the same operators present in the model formulation.
For the OSRC preconditioner, the hyperparameters are given by a Padé series with four terms, a branch cut of angle $\pi/3$, sparse $LU$ decompositions, and a damping parameter of $\epsilon = 0.4(k_m R_m)^{-2/3}$. Different values are used for the damping, including the interior and exterior wavenumbers.

\paragraph{Software and hardware platform}

The BEM was implemented with version~3.3 of the open-source BEMPP library~\cite{smigaj2015solving, scroggs2017software}. The library SciPy version 1.2.1 was used for the linear algebra~\cite{scipy}. Graphics were created with the libraries Matplotlib~\cite{hunter2007matplotlib} and Seaborn~\cite{waskom2020seaborn}.
All simulations were performed with hyperthreading on a workstation with 16 processor cores (Intel\textregistered~Xeon(R) CPU E5-2683 v4\textcircled{a}2.10 GHz) and 512 GB RAM. Shared-memory parallelisation was achieved through threaded Lapack routines called by SciPy and an Intel TBB implementation in BEMPP~\cite{smigaj2015solving}.
Additional gains in compute time can be achieved through high-performance computing on graphics cards~\cite{molina2018iterative} or clusters~\cite{yokota2012tuned}, which is outside the scope of this study.

\subsection{Performance studies}

The benchmarking includes a total number of 538 different preconditioned boundary integral formulations and a variety of different geometries. Below, computational results of the benchmarks will be presented.

\subsubsection{Matrix assembly}

\begin{figure}[!t]
	\centering
	\includegraphics[width=\columnwidth]{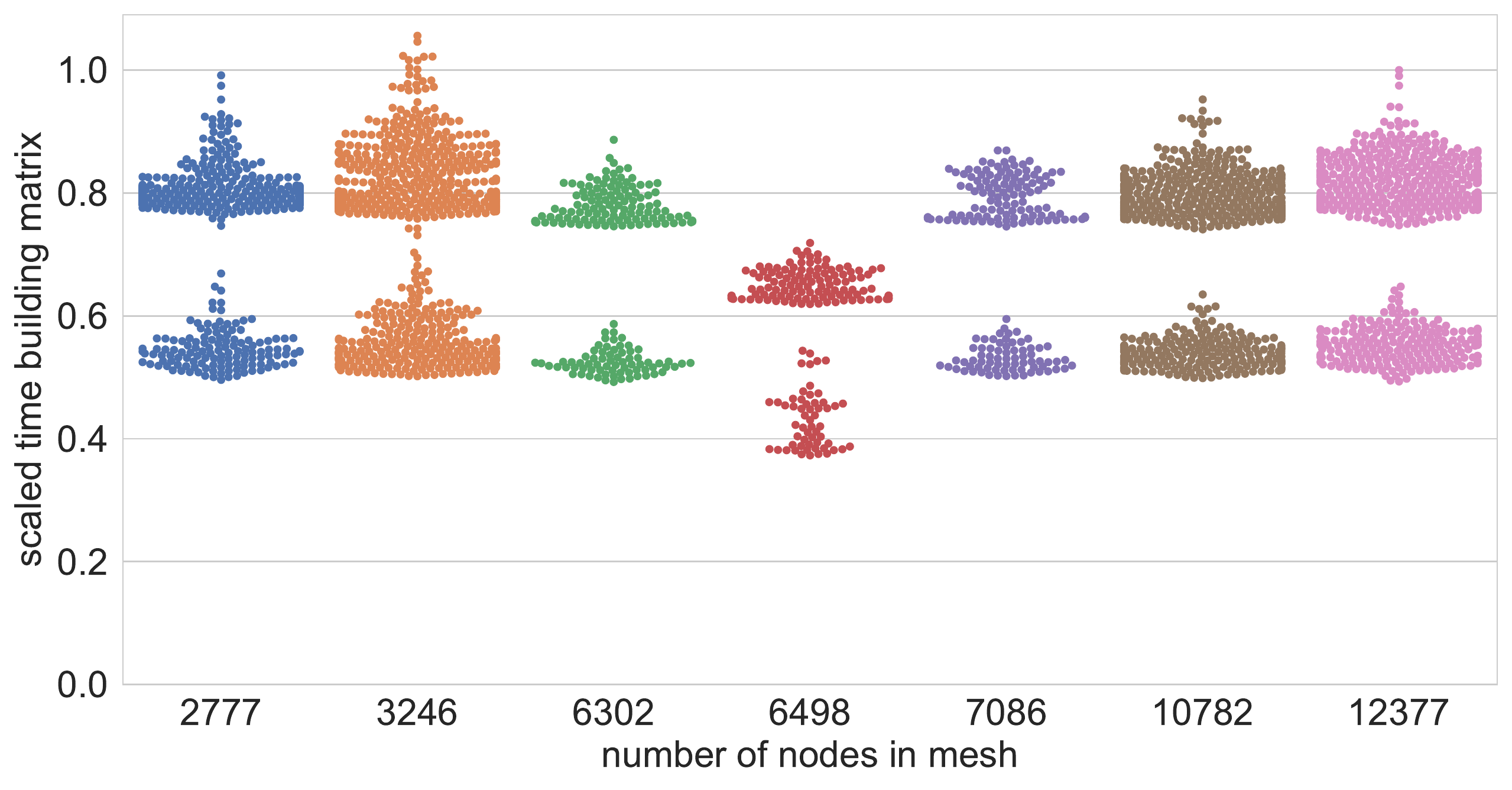}
	\caption{A swarm plot of the timing of the matrix assembly, scaled with the square of the number of nodes in the mesh, and normalised by the maximum time at 12\,377~nodes.}
	\label{fig:matrixassembly}
\end{figure}

Let us first benchmark the time to assemble the model. For dense matrices, the computational complexity is $\mathcal{O}(n^2)$ with $n$ the number of degrees of freedom, which yields a scaling of $\mathcal{O}(f^4)$ for a fixed number of elements per wavelength. This computational complexity is confirmed by the benchmarks, where the measured time to build the linear system has a scaling of 1.98 with respect to the number of degrees of freedom. The benchmark results presented in Figure~\ref{fig:matrixassembly} clearly show two groups of formulations at each mesh. This observation is consistent with the number of boundary integral operators in the formulation, which is either two or four, as summarised in Table~\ref{table:formulations}. Notice that a ratio of two to three is observed since the adjoint double-layer operator is not explicitly assembled: the transpose of the double-layer operator was used. Furthermore, the sudden drop visible at 6498~nodes is expected because this benchmark corresponds to the case of two spheres and fewer operators are necessary for the cross interactions.

\FloatBarrier
\subsubsection{Time per GMRES iteration}

\begin{figure}[!t]
	\centering
	\includegraphics[width=\columnwidth]{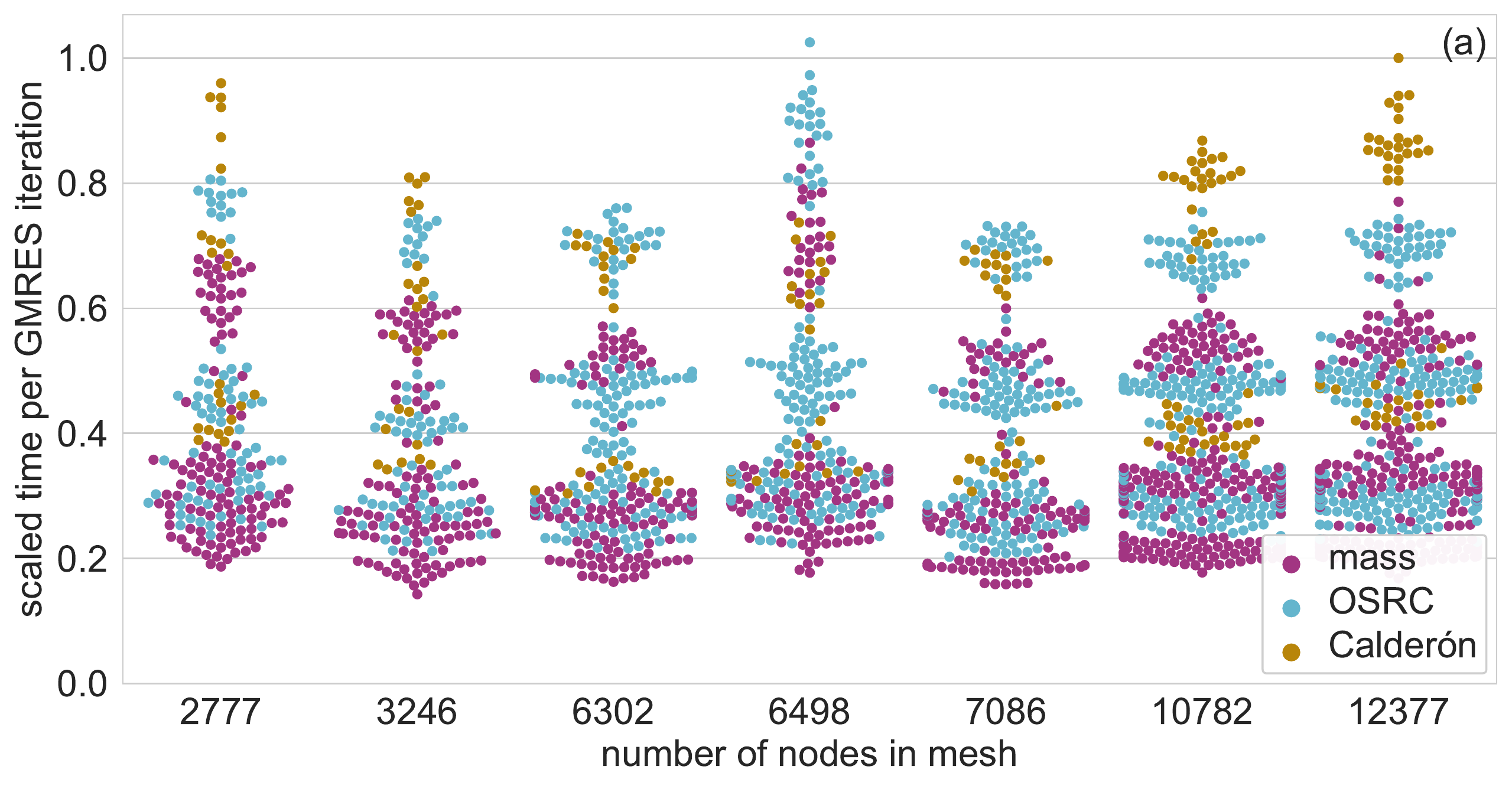}
	\includegraphics[width=\columnwidth]{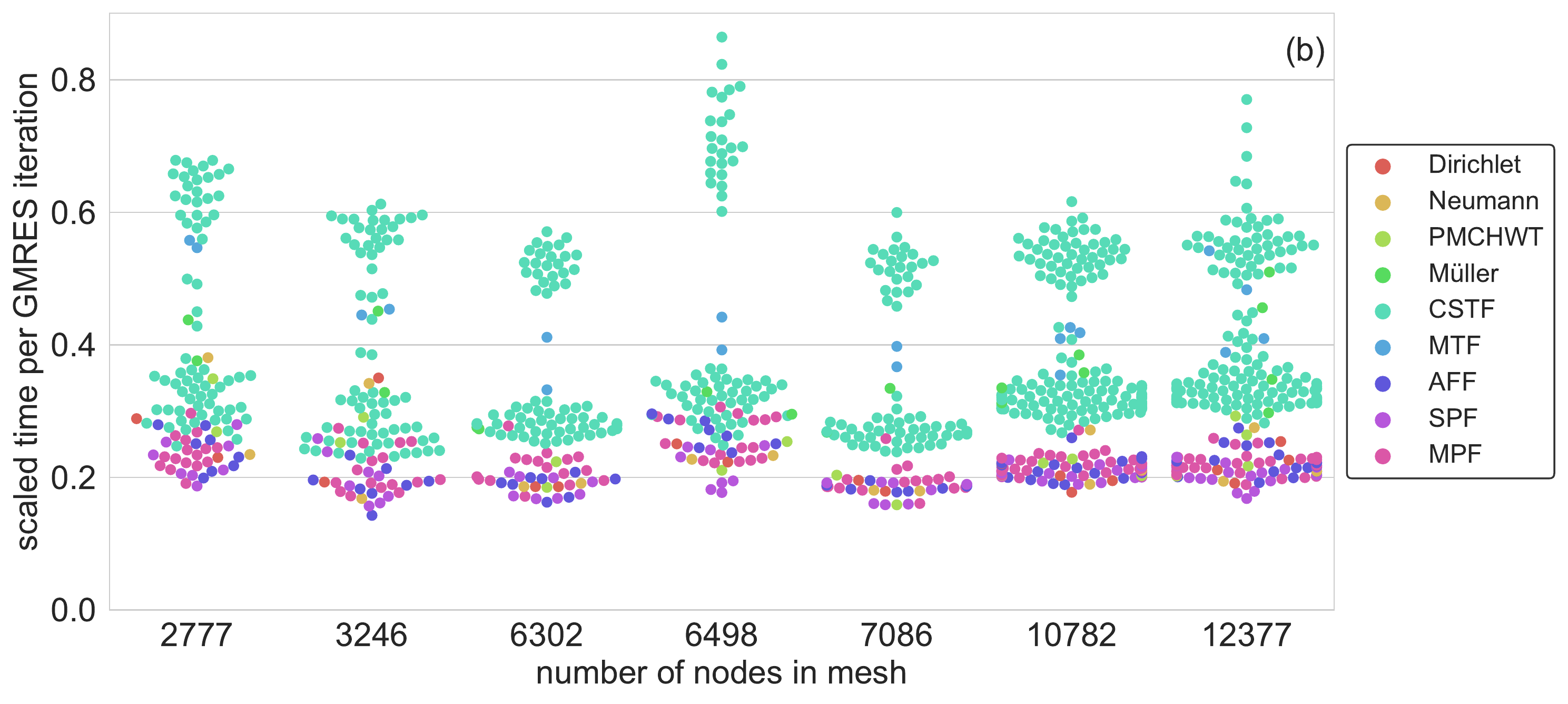}
	\caption{A swarm plot of the time per GMRES iteration, averaged over the GMRES history of each simulation, scaled with $N_\text{nodes}^\alpha$ for $\alpha=1.28$, and normalised by the maximum time at 12\,377~nodes. (a) Separated by preconditioner type. (b) Separated by formulation for mass preconditioning only.}
	\label{fig:matvec}
\end{figure}

Since the BEM requires dense matrix arithmetic, the time per GMRES iteration is dominated by the multiplication of the preconditioned matrix with a vector. No preconditioner system needs to be solved because the Calderón preconditioning is a matrix-vector multiplication and the OSRC preconditioner uses a set of sparse LU decompositions calculated at the assembly stage.
The expected scaling of a dense matrix-vector multiplication is $\mathcal{O}(n^2)$. Surprisingly, the benchmark suggests a scaling of $N_\text{nodes}^\alpha$ for $\alpha=1.28$, as can be seen in Figure~\ref{fig:matvec}. The observed scaling is much better than expected, which is likely because of the optimised Lapack routines for linear algebra routines~\cite{golub2013matrix}. This is in contrast to the matrix assembly, which is performed by special quadrature rules implemented in the C++ kernel of the BEMPP library~\cite{smigaj2015solving}.
These timing characteristics strongly depend on the design choices of the software package, and will also drastically change when accelerators like the fast-multipole method or hierarchical matrix compression are employed.

When comparing the different preconditioned formulations, Figure~\ref{fig:matvec}(a) confirms the expected increase in time per GMRES iteration with preconditioning, with Calderón preconditioning the most expensive one due to its dense blocks. As before, the groups of timings correspond to formulations with the same number of operators involved, which is also visible in Figure~\ref{fig:matvec}(b). Notice that although the results in Figure~\ref{fig:matvec}(b) are for the mass preconditioner only, this still includes combined single-trace formulation with the OSRC operator as coupling operator, which are the most expensive cases.

\FloatBarrier
\subsubsection{Accuracy of field reconstruction}
\label{sec:benchmark:error}

For transmission through a single spherical object, the analytical solution is given by a series expansions in spherical harmonics. With the purpose of an accuracy analysis, the total pressure field~$\utot$ is evaluated on a visualisation grid of $101 \times 101$ points that are uniformly located on a square of size $3 \times 3$\,cm, centered at the origin of the sphere and with nodes both in the exterior and interior of the sphere. The accuracy is defined by the peak signal-to-noise ratio (PSNR), which is a common quality measure for image reconstruction and defined as the ratio of the mean squared error and the maximum, in decibels:
\begin{equation}
	\mathrm{PSNR}
	= -10 \cdot \log_{10} \left(\frac{\frac1N \sum_i \left( \utot^\mathrm{exact}(\mathbf{x}_i) - \utot(\mathbf{x}_i) \right)^2}{\left(\max_i \utot^\mathrm{exact}(\mathbf{x}_i)\right)^2}\right)
\end{equation}
where $\mathbf{x}_i$ for $i=1,2,\dots,N$ are the field points and $\utot^\mathrm{exact}$ denotes the analytical solution.

\begin{figure}[!t]
	\centering
	\includegraphics[width=\columnwidth]{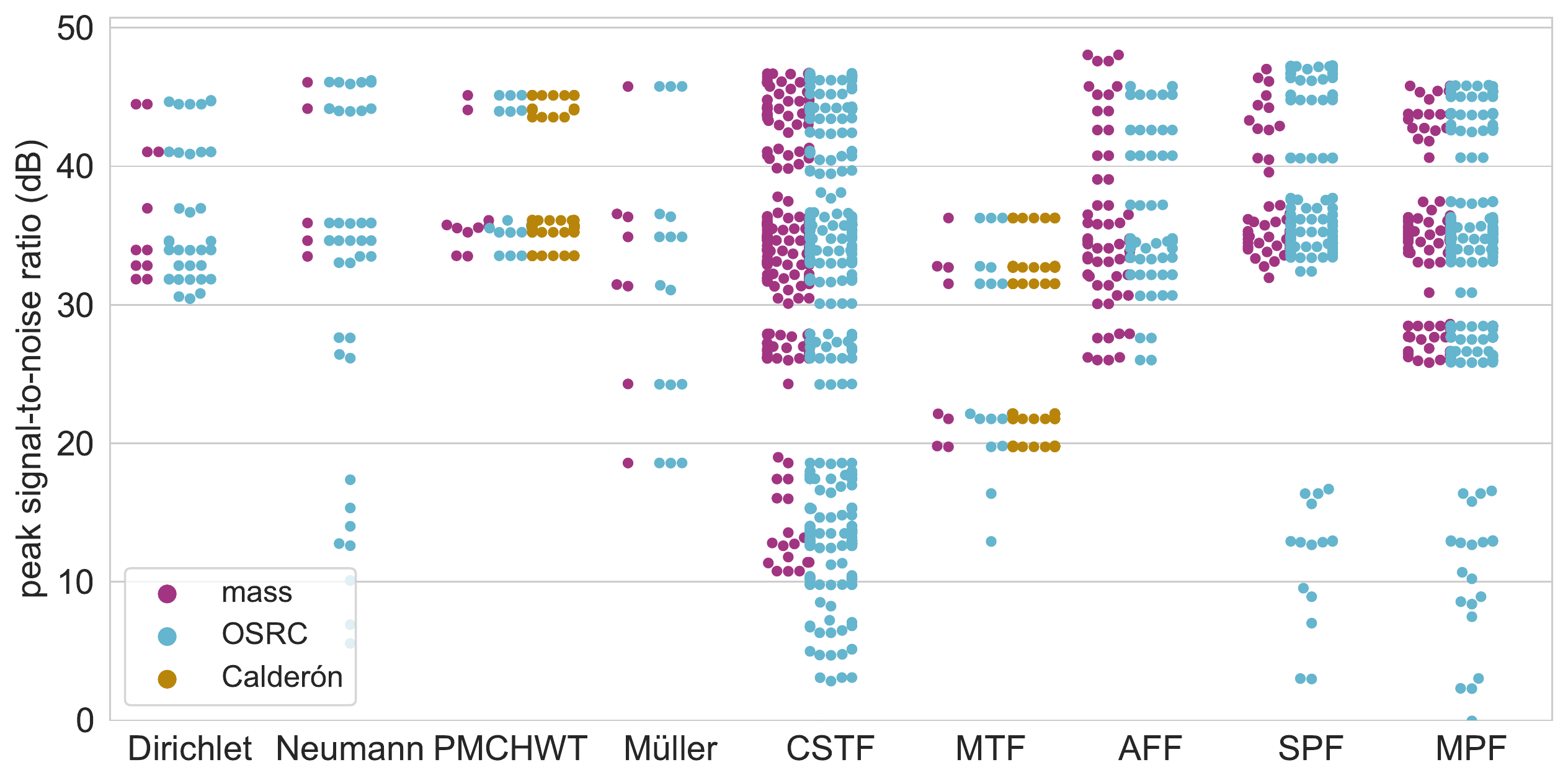}
	\caption{A swarmplot of the accuracy for all preconditioned formulations on one sphere that converged within 1000 iterations. Higher values of PSNR mean higher accuracy of the field reconstruction. As a reference, a PSNR accuracy of 20~dB and 40~dB correspond to L2 errors of around 5\% and 1\%, respectively.}
	\label{fig:error}
\end{figure}

The reconstruction accuracy measured in the benckmark is presented in Figure~\ref{fig:error}, which shows large differences between the preconditioned formulations, even though all simulations use the same GMRES tolerance. The GMRES implementation of the benchmarks take the preconditioned matrix norm of the residual as termination criterion. This is a different error measure than the PSNR of the pressure fields. Overall, the results show that preconditioning does not change the quality of the fields, except for several cases with poorly designed preconditioned formulations. In these special cases, the fields are inaccurate even when GMRES converges. Simulations with a PSNR lower than 20~dB have visually appreciable deficiencies in the fields and will be excluded from analysis in the following sections. Furthermore, the accuracy depends on the material parameters and increasing the mesh density improves the accuracy.

\FloatBarrier
\subsubsection{Convergence with material parameters and frequency}

\begin{figure}[!t]
	\centering
	\includegraphics[width=\columnwidth]{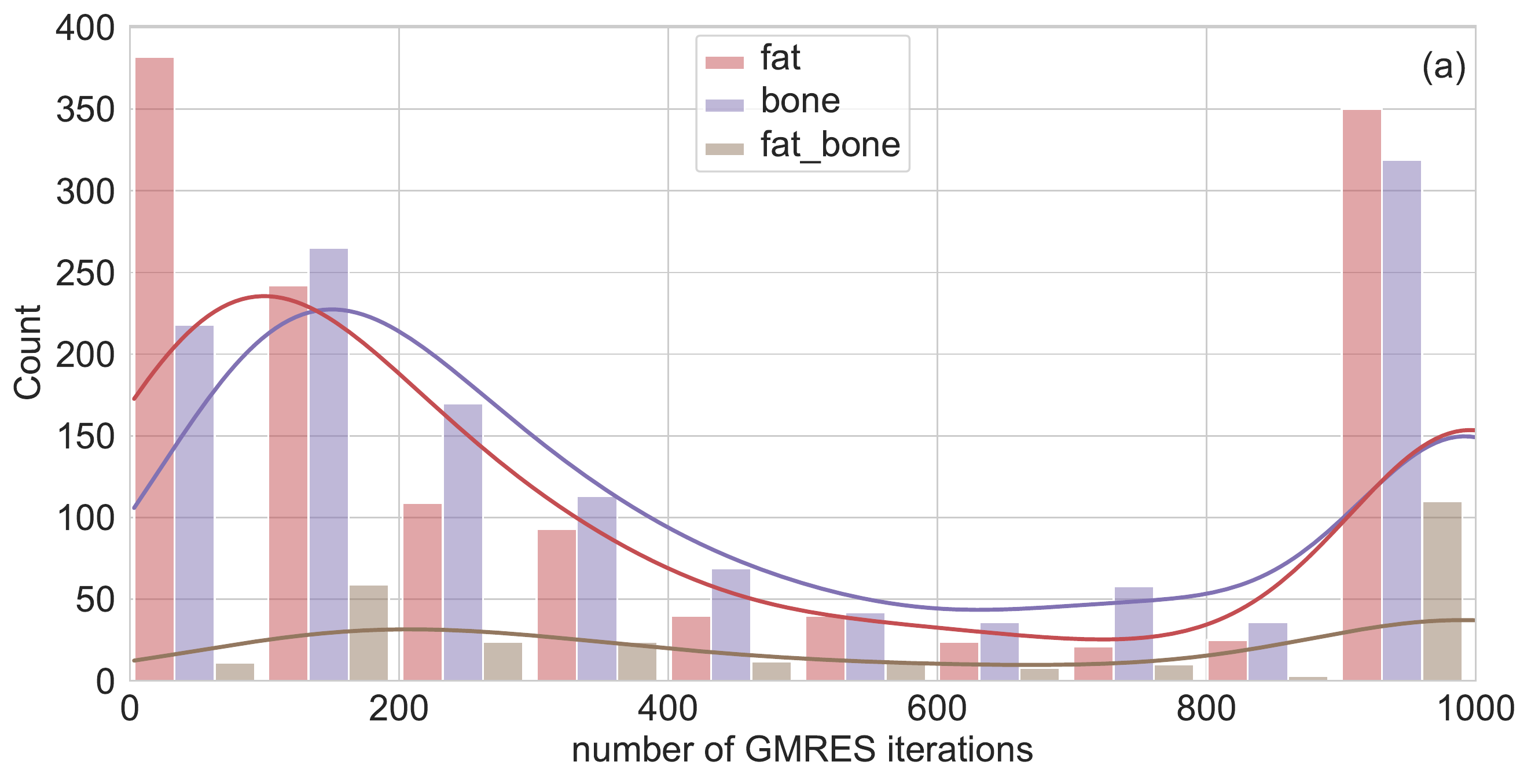}
	\includegraphics[width=\columnwidth]{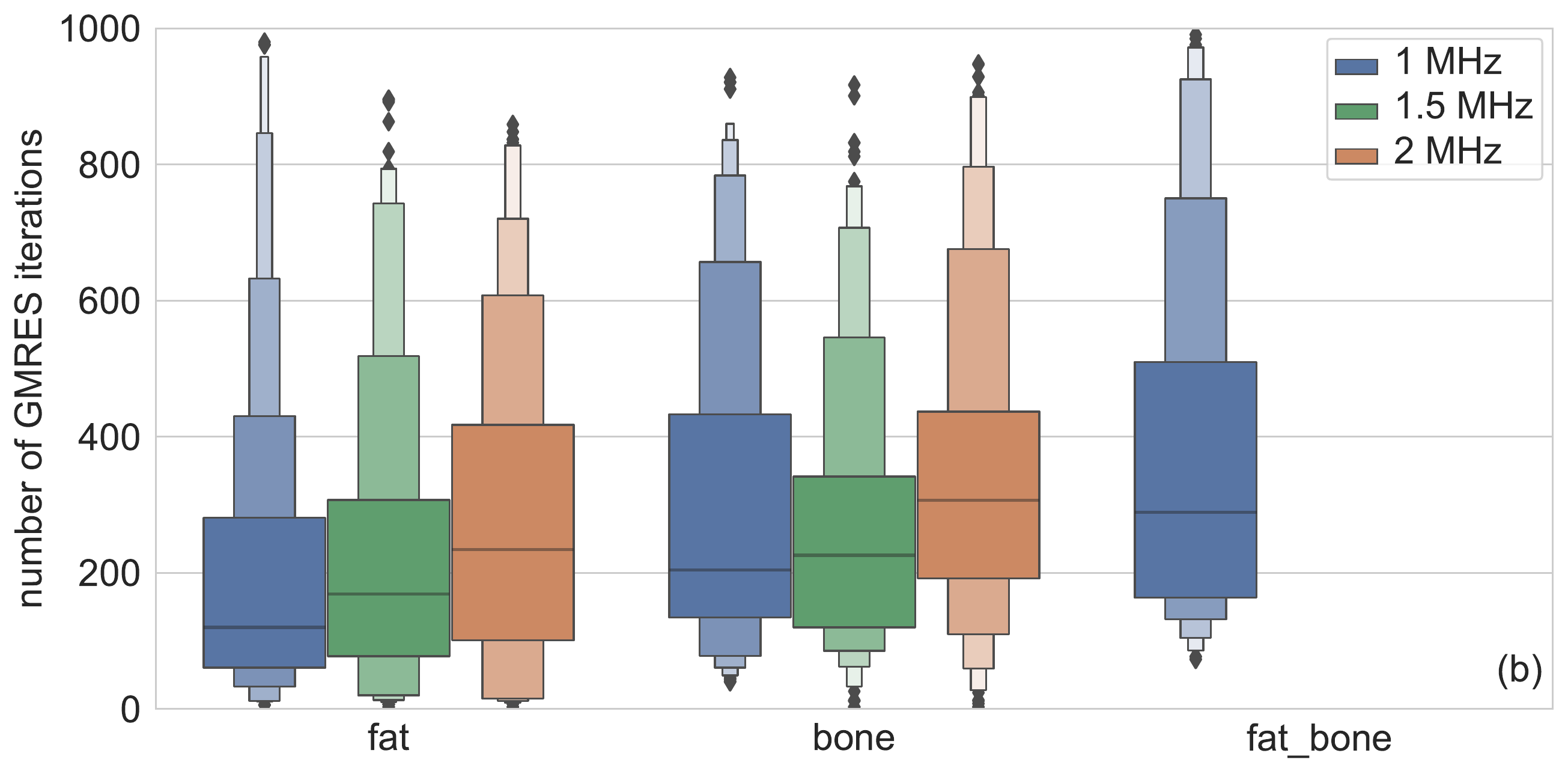}
	\caption{The number of GMRES iterations for each preconditioned formulation, with the maximum set to 1000. (a) A histogram and kernel density estimation for all formulations. The count is the number of preconditioned formulations inside each bin that consists of an interval of 100 GMRES iterations. The last bin includes formulations that did not converge in 1000 GMRES iterations. (b) A boxenplot (or letter-value plot~\cite{heike2017letter}) for the formulations that converged. From the median at the horizontal line, the consecutive blocks (upwards and downwards) contain roughly 25\%, 12.5\%, 6.25\%, etc.~of the preconditioned formulations. Outliers are presented by diamonds.}
	\label{fig:gmres:material}
\end{figure}

The convergence of GMRES strongly depends on the material characteristics and the driving frequency, as confirmed by the benchmark results presented in Figure~\ref{fig:gmres:material}. The large proportion at the right in Figure~\ref{fig:gmres:material}(a) includes all simulations that did not yet converge after the maximum of 1000~iterations. Figure~\ref{fig:gmres:material}(b) only includes the benchmark results for which GMRES converged within 1000 iterations. These results show a skewed distribution of the number of GMRES iterations for the preconditioned formulations. That is, the median is low but a significant proportion of preconditioned formulations requires considerably more iterations for GMRES to converge.

The dependency of the GMRES convergence on the material characteristics can be attributed to different mechanisms. For example, bone has more attenuation, a longer wavelength, and requires less degrees of freedom compared to fat, which can all be considered as favourable. However, fat has a smaller contrast with water in terms of wavespeed and density, as compared to bone and water. The benchmarking suggests that the contrasts in material parameters are the dominant characteristics for the convergence of GMRES. The observed deterioration of convergence when frequency increases is expected and specialised formulations and preconditioners need to be designed to improve the convergence at high frequencies~\cite{betcke2017computationally}.

The results presented in Figure~\ref{fig:gmres:material} do not distinguish between preconditioned formulations and, therefore, any conclusion drawn from these benchmarks depend on the stratification of formulation parameters that were chosen. As will be shown below, these results are general statements that do not necessarily hold for a specific preconditioned formulation. In any case, these benchmarks give important information when preconditioned formulations are chosen without optimising the performance for the specific configuration.

\FloatBarrier
\subsubsection{Convergence of formulations}

\begin{figure}[!t]
	\centering
	\includegraphics[width=\columnwidth]{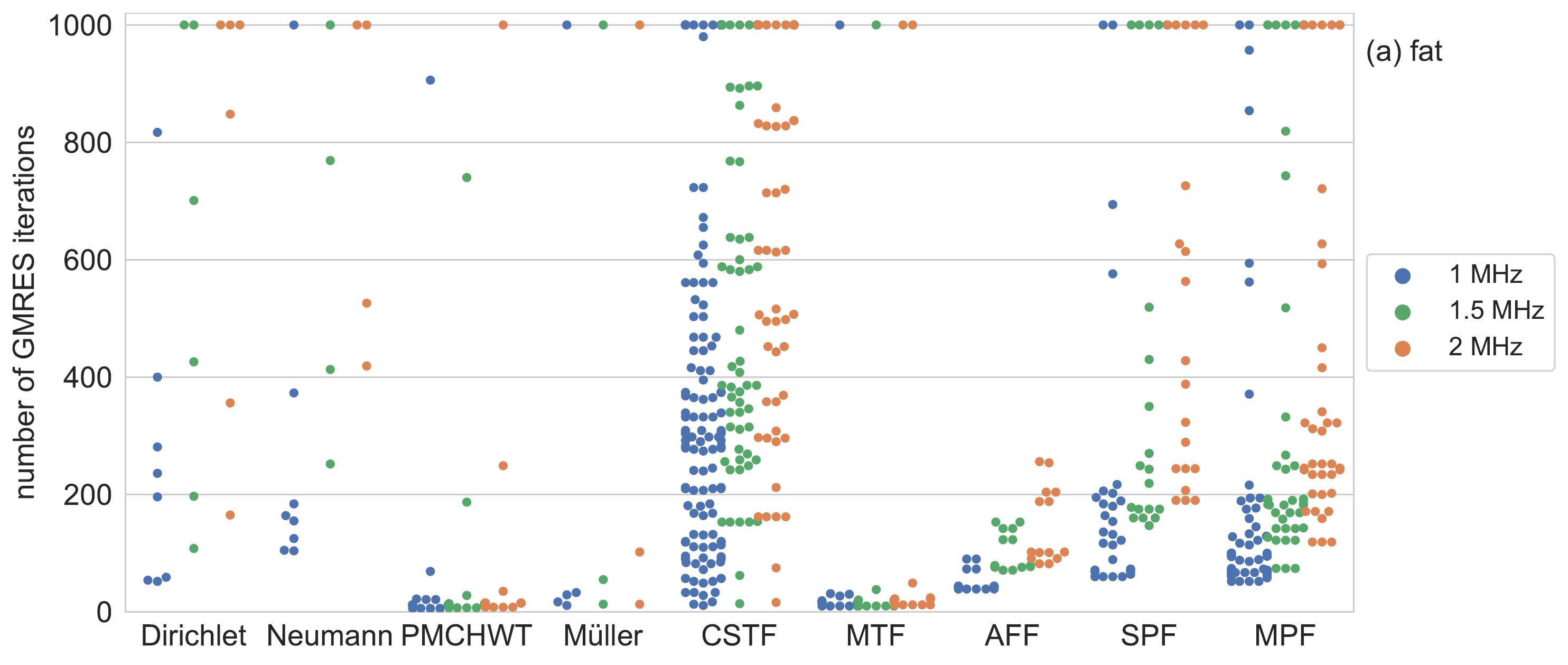}
	\includegraphics[width=\columnwidth]{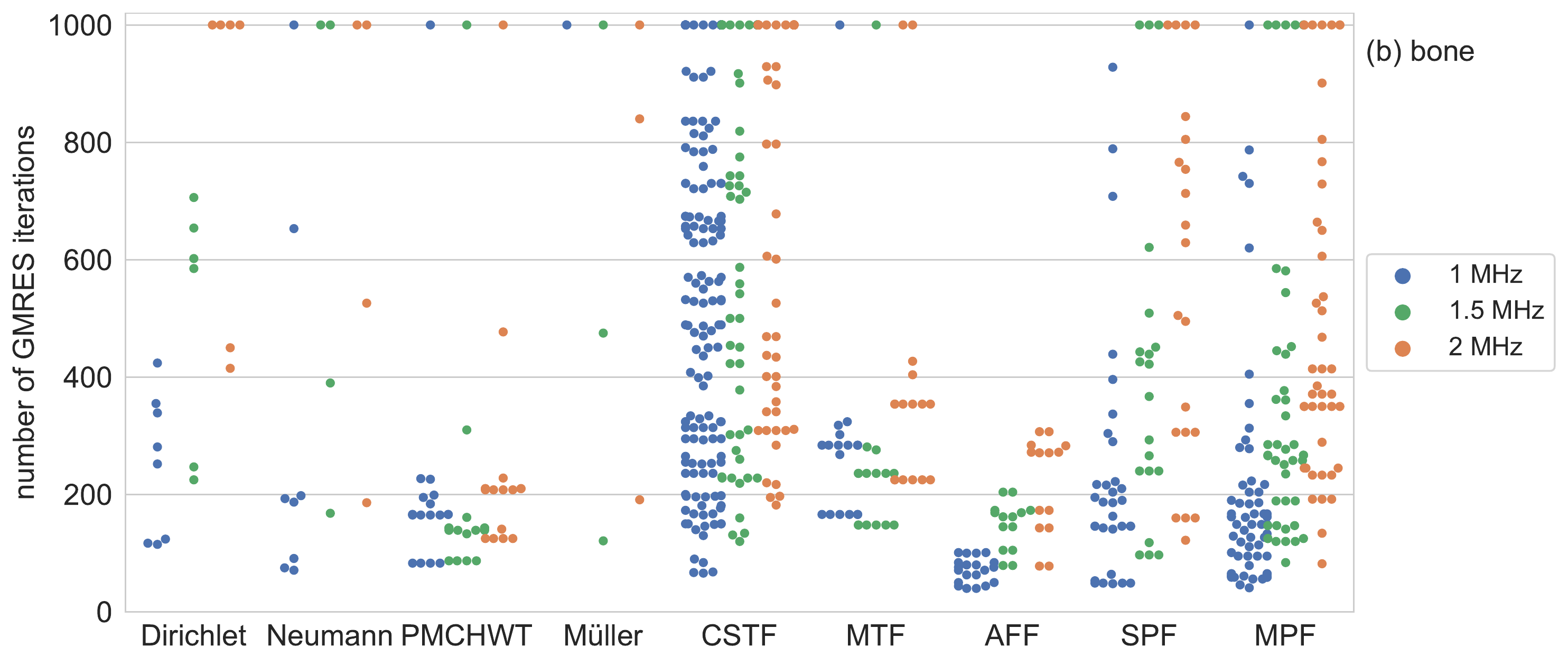}
	\caption{A swarm plot of the number of GMRES iterations for each preconditioned formulation, with the maximum set to 1000. The material is (a) fat and (b) bone.}
	\label{fig:gmres:model}
\end{figure}

The benchmarking results presented in Figure~\ref{fig:gmres:model} demonstrate that the convergence of the BEM strongly depends on the boundary integral formulation and preconditioner. There is a multifaceted interaction between the frequency, material type and formulation. For example, the formulations AFF, SPF and MPF have similar convergence behaviour for both material types, whereas the PMCHWT and MTF have a sharp increase in number of iterations. Furthermore, a wide spread in convergence behaviour is observed with some formulations converging within a few iterations while others did not converge after a thousand iterations. This confirms that choosing the correct boundary integral formulation is essential to obtain efficient simulations.

\FloatBarrier
\subsubsection{Convergence of preconditioners}

\begin{figure}[!t]
	\centering
	\includegraphics[width=\columnwidth]{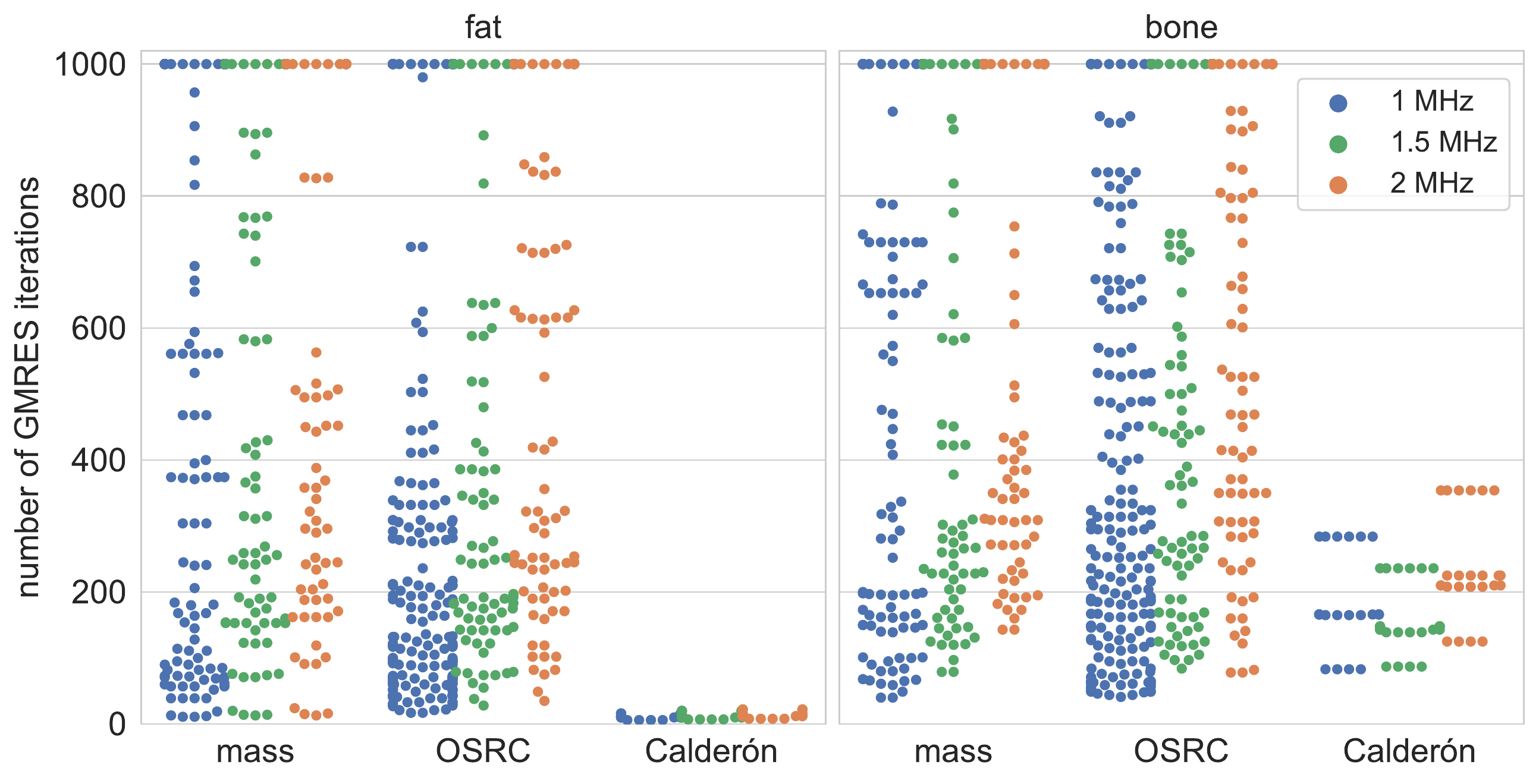}
	\caption{A swarmplot of the number of GMRES iterations for each preconditioned formulation on a single sphere, with the maximum set to 1000.}
	\label{fig:gmres:preconditioner}
\end{figure}

The results in Figure~\ref{fig:gmres:preconditioner} confirm that preconditioning works in general, since most of the OSRC and Calderón preconditioned formulations require less GMRES iterations than the simple mass preconditioner. The few cases where preconditioning deteriorates the convergence correspond to OSRC preconditioned formulations that use poorly chosen hyperparameters. The Calderón preconditioning is robust and requires few iterations, especially for fat. However, the convergence behaviour of Calderón preconditioning should consider the multifaceted character of preconditioned boundary integral formulations: Calderón preconditioner is available for the PMCHWT and MTF formulations only. Both formulations are stable but with convergence issues when materials have a high contrast across the interface. Also, preconditioning has an influence on the time per iteration, as was already discussed in Figure~\ref{fig:matvec}.

\FloatBarrier
\subsubsection{Multiple objects}

\begin{figure}[!t]
	\centering
	\includegraphics[width=\columnwidth]{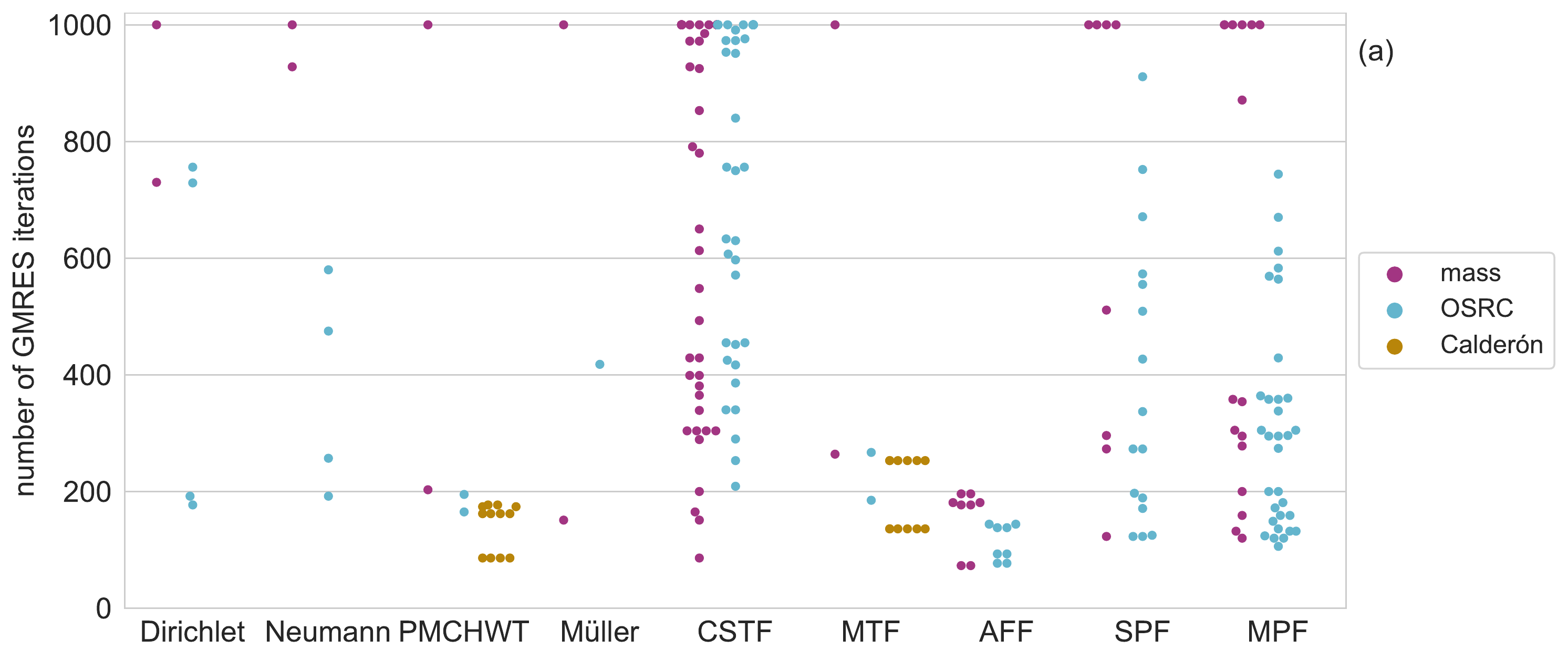}
	\includegraphics[width=\columnwidth]{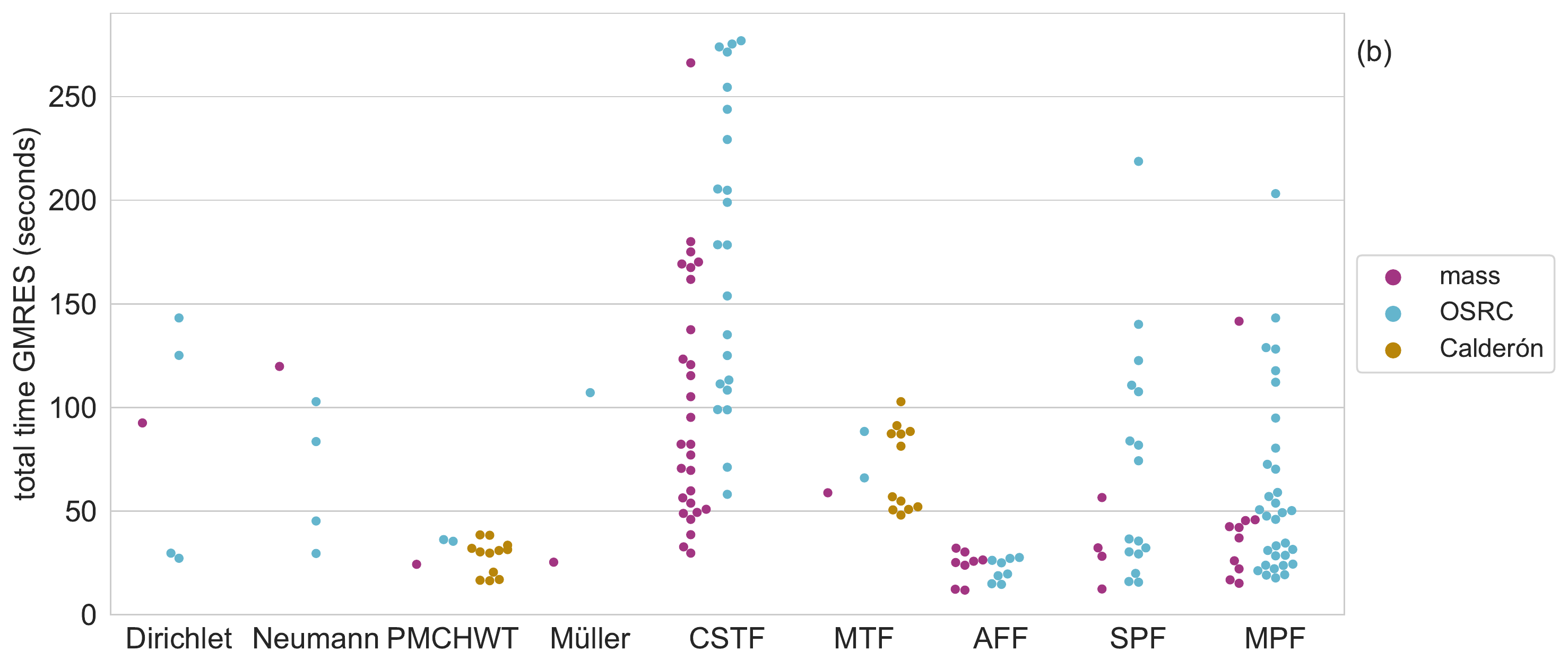}
	\caption{A swarmplot of (a) the number of iterations and (b) the computation time of GMRES, for the preconditioned formulations on two spheres. In (b) only formulations that converged are depicted.}
	\label{fig:gmres:twosphere}
\end{figure}

The results in Figure~\ref{fig:gmres:twosphere} present the GMRES performance in the case of the two-sphere benchmark test. Again, the convergence can improve considerably with preconditioning, but only when correctly designed. The benchmarks also show the trade-off between less iterations and more computation time per iterations with preconditioning. Considering the total time to solve the linear system, the indirect formulations are relatively efficient since they involve less boundary integral operators.

\FloatBarrier
\subsection{Large-scale simulation}

\begin{figure}[!t]
	\centering
	\includegraphics[width=\columnwidth]{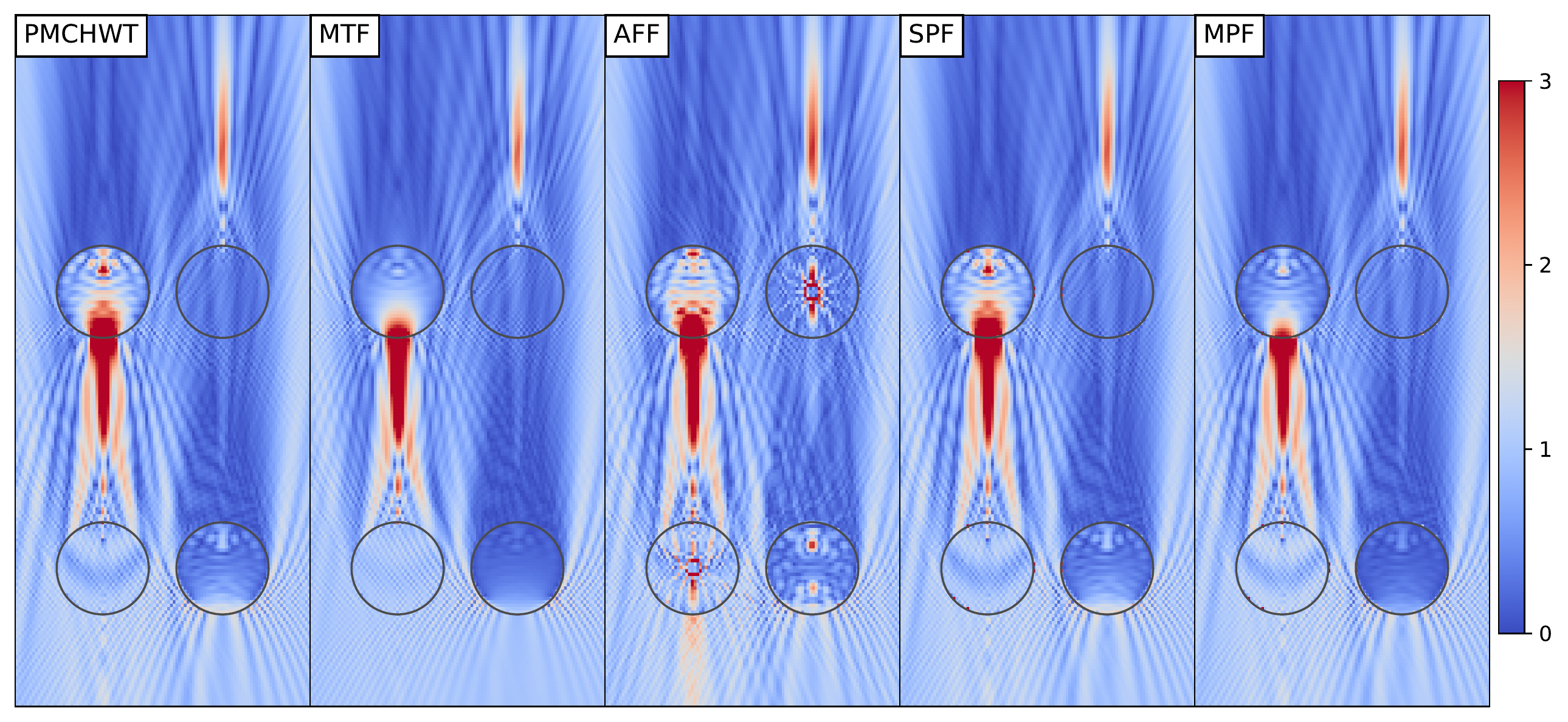}
	\caption{The amplitude of the acoustic field on the plane given by $z=0$. The incident plane wave field travels upwards. The lower left and upper right spheres are of fat material and the upper left and lower right spheres are bone.}
	\label{fig:fourspheres:field}
\end{figure}

The previous benchmarks were performed on test cases with an intermediate complexity feasible to test hundreds of different formulations. Now, let us consider several of the best performing formulations and test them on a large-scale geometry. Four spheres with a radius of 5~mm are embedded in an exterior medium of water, two of the spheres are bone and the other two fat, with centers at $(\pm 6.5, \pm 15, 0)$~mm. The incident plane wave field travels in positive $y$-direction and has a driving frequency of 2~MHz. The five elements per wavelength result in a mesh of 72\,254 nodes. In order to fit the matrices in the memory, hierarchical matrix compression was used with a tolerance of $10^{-6}$. For the OSRC preconditioner, the standard parameter settings were used and for the Calderón preconditioner, the full version was used. The field is visualised in Figure~\ref{fig:fourspheres:field}. As expected, the spheres made of fat, which is a soft material, create a lensing effect whereas the spheres made of bone, which is a hard material, create a shadow region.

\begin{table}[!t]
	\caption{The performance characteristics of the large-scale benchmark. The time ($T$) is divided into building and solving the linear system.}
	\label{table:fourspheres:performance}
	\centering
	\begin{tabular}{llrrrrr}
		\hline\hline
		model & preconditioner 
		& $T$ build & \#iter & $T$ solve & $T/$iter \\
		\hline
		PMCHWT & Calderón
		& 6:51 h &  281 & 0:55 h & 11.77 s \\
		PMCHWT & OSRC
		& 6:58 h &  224 & 0:24 h & 6.56 s \\
		MTF & Calderón
		& 6:45 h &  579 & 2:03 h & 12.70 s \\
		AFF \eqref{eq:formulation:auxiliary:neu:phi} & OSRC
		& 5:41 h &  285 & 0:24 h & 5.01 s \\
		SPF \eqref{eq:formulation:singlepotential:dlpext:dlpint} & OSRC
		& 4:40 h & 1000 & 1:14 h & 4.46 s \\
		MPF \eqref{eq:formulation:mixedpotential:dlp:ntd} & OSRC
		& 5:16 h &  427 & 0:35 h & 4.86 s \\
		\hline\hline
	\end{tabular}
\end{table}

Table~\ref{table:fourspheres:performance} summarises the performance statistics for the large-scale benchmark. The PMCHWT and MTF both require long build times since they use all boundary integral operators. The other formulations use less operators and are quicker in the matrix assembly, where the timing of the hierarchical matrix compression  also depends on the specific operators present in the formulation. The additional time for building the OSRC preconditioner is less than 30~seconds in all cases whereas Calderón preconditioning does not require additional assembly time. However, Calderón preconditioning doubles the time per iteration whereas the sparse OSRC preconditioner has little overhead.

The SPF is the only model that did not converge in 1000 iterations while the PMCHWT and AFF have the smallest numbers of iterations. Notice that the AFF is the quickest in overall time but, as can be seen in Figure~\ref{fig:fourspheres:field}, the field is inaccurate. Even though GMRES converged, spurious solutions are present in the interior of the spheres. Performing the same model with a higher mesh resolution solved this issue, of course at the expense of considerably longer computation times. A similar behaviour was observed in Section~\ref{sec:benchmark:error} as well: a small error in the matrix norm for GMRES does not necessarily result in small errors in the pressure field. Differently, while the SPF did not converge, the pressure field is accurately retrieved. The PMCHWT formulation is very robust and simulates acoustic fields accurately, while the OSRC preconditioner yields fast convergence at high frequencies.

\section{Conclusions}

This study surveyed the design of five families of boundary integral formulations for acoustic transmission through multiple domains. Each of these formulations use a different potential representation of the field or a different coupling algorithm at the material interface. Calderón and OSRC preconditioning were applied to all feasible formulations, leading to novel combinations of operator preconditioning and boundary integral formulations. Extensive benchmarks compare the computational performance of hundreds of preconditioned boundary integral formulations in terms of solution accuracy, GMRES convergence, and calculation time. The numerical results confirm that a proper choice and correct design of the model can improve the calculation time to solve the discretised system by orders of magnitude. Furthermore, there is a multifaceted dependency of the performance on material type, driving frequency and multiple reflection. None of the preconditioned formulations outperforms all others on each benchmark considered. Instead, the methodology needs to be adjusted to the specific configuration of the model, such as frequency and material types, as well as the efficiency objectives, such as memory consumption, calculation time, accuracy and robustness. Hence, expert knowledge is required, and this study presented the major considerations to be taken into account by a BEM practitioner.

The benchmarks show general recommendations on the choice of preconditioned boundary integral formulations. Figures~\ref{fig:matrixassembly} and~\ref{fig:matvec} show that dense matrix arithmetic is quick for small-scale problems, but compression techniques are required when more than ten thousand nodes are present in the surface mesh. Figure~\ref{fig:error} shows that the PMCHWT is one of the most robust formulations, with highly accurate field reconstructions for only four elements per wavelength and a moderate GMRES tolerance. Figures~\ref{fig:gmres:material} and~\ref{fig:gmres:model} show that the indirect formulations converge quickly for high-contrast materials. Figures~\ref{fig:gmres:preconditioner} and~\ref{fig:gmres:twosphere} show that mass-matrix preconditioning is sufficient at low frequencies, and OSRC preconditioning is very effective at high frequencies and multiple reflection. Finally, the large-scale simulations in Figure~\ref{fig:fourspheres:field} show that indirect formulations have short assembly time but slow convergence and inaccuracies, and the OSRC-preconditioned PMCHWT formulation is efficient and robust.

As with any benchmarking study, the parameter space was restricted due to practical limitations. Interesting dependencies on the geometry that were not considered for brevity include nonsmooth domains and resonant cavities. Furthermore, this study focuses only on the Helmholtz equation for acoustic wave propagation. Similar strategies can design formulations for electromagnetics and elastodynamics. Finally, even though the number of preconditioned formulations considered in this study is impressive, it is by no means exhaustive. The design freedom for boundary integral formulations and preconditioners allows for the development of highly specialised techniques that are optimised for a specific setting. These might outperform the preconditioned formulations presented here. In general, these benchmarks of preconditioned formulations will suffice for the creation of robust and efficient boundary element methods for most practical purposes.

\section*{Acknowledgment}

This work was financially supported by CONICYT [FONDECYT 11160462], the Vicerrectoría de Investigación of the Pontificia Universidad Católica de Chile, and the EPSRC [EP/P012434/1].

\bibliographystyle{unsrt}
\bibliography{refs}

\begin{thebibliography}{10}

\bibitem{lahaye2017modern}
D.~Lahaye, J.~Tang, and K.~Vuik.
\newblock {\em Modern Solvers for {H}elmholtz Problems}.
\newblock Birkh{\"a}user, Cham, 2017.

\bibitem{nedelec2001acoustic}
Jean-Claude N{\'e}d{\'e}lec.
\newblock {\em Acoustic and electromagnetic equations: integral representations
  for harmonic problems}.
\newblock Springer, New York, 2001.

\bibitem{steinbach2008numerical}
Olaf Steinbach.
\newblock {\em Numerical approximation methods for elliptic boundary value
  problems: finite and boundary elements}.
\newblock Springer, New York, 2008.

\bibitem{hsiao2008boundary}
George~C Hsiao and Wolfgang~L Wendland.
\newblock {\em Boundary integral equations}.
\newblock Springer, Berlin, 2008.

\bibitem{sauter2010boundary}
Stefan~A Sauter and Christoph Schwab.
\newblock {\em Boundary Element Methods}.
\newblock Springer, Berlin, 2010.

\bibitem{greengard1987fast}
Leslie Greengard and Vladimir Rokhlin.
\newblock A fast algorithm for particle simulations.
\newblock {\em Journal of Computational Physics}, 73(2):325--348, 1987.

\bibitem{hackbusch2015hierarchical}
Wolfgang Hackbusch.
\newblock {\em Hierarchical matrices: algorithms and analysis}.
\newblock Springer, Berlin, 2015.

\bibitem{baydoun2018quantification}
Suhaib~Koji Baydoun and Steffen Marburg.
\newblock Quantification of numerical damping in the acoustic boundary element
  method for two-dimensional duct problems.
\newblock {\em Journal of Theoretical and Computational Acoustics},
  26(03):1850022, 2018.

\bibitem{smigaj2015solving}
Wojciech {\'S}migaj, Timo Betcke, Simon Arridge, Joel Phillips, and Martin
  Schweiger.
\newblock Solving boundary integral problems with {BEM++}.
\newblock {\em ACM Transactions on Mathematical Software (TOMS)}, 41(2):6,
  2015.

\bibitem{chew2001fast}
Weng~Cho Chew, Eric Michielssen, JM~Song, and Jian-Ming Jin.
\newblock {\em Fast and efficient algorithms in computational
  electromagnetics}.
\newblock Artech House, Inc., Norwood, MA, 2001.

\bibitem{marburg2018boundary}
Steffen Marburg.
\newblock Boundary element method for time-harmonic acoustic problems.
\newblock In {\em Computational Acoustics}, pages 69--158. Springer, 2018.

\bibitem{kirkup2019boundary}
Stephen Kirkup.
\newblock The boundary element method in acoustics: A survey.
\newblock {\em Applied Sciences}, 9(8):1642, 2019.

\bibitem{maue1949formulierung}
A.-W. Maue.
\newblock Zur {F}ormulierung eines allgemeinen {B}eugungs-problems durch eine
  {I}ntegralgleichung.
\newblock {\em Zeitschrift f{\"u}r Physik}, 126(7-9):601--618, 1949.

\bibitem{muller1957grundprobleme}
Claus M{\"u}ller.
\newblock {\em Grundprobleme der mathematischen {T}heorie elektromagnetischer
  {S}chwingungen}.
\newblock Springer, Berlin, 1957.

\bibitem{mitzner1966acoustic}
Kenneth~M Mitzner.
\newblock Acoustic scattering from an interface between media of greatly
  different density.
\newblock {\em Journal of Mathematical Physics}, 7(11):2053--2060, 1966.

\bibitem{mautz1977electromagnetic}
Joseph~R Mautz and Roger~F Harrington.
\newblock Electromagnetic scattering from a homogeneous body of revolution.
\newblock Technical report, Syracuse University, Syracuse, NY, 1977.
\newblock Technical Report TR-77-10.

\bibitem{costabel1985direct}
Martin Costabel and Ernst Stephan.
\newblock A direct boundary integral equation method for transmission problems.
\newblock {\em Journal of Mathematical Analysis and Applications},
  106(2):367--413, 1985.

\bibitem{acosta2015surface}
Sebastian Acosta.
\newblock On-surface radiation condition for multiple scattering of waves.
\newblock {\em Computer Methods in Applied Mechanics and Engineering},
  283:1296--1309, 2015.

\bibitem{alzubaidi2016formulation}
Hasan Alzubaidi, Xavier Antoine, and Chokri Chniti.
\newblock Formulation and accuracy of on-surface radiation conditions for
  acoustic multiple scattering problems.
\newblock {\em Applied Mathematics and Computation}, 277:82--100, 2016.

\bibitem{poggio1973integral}
A.~J. Poggio and E.~K. Miller.
\newblock Integral equation solutions of three-dimensional scattering problems.
\newblock In R.~Mittra, editor, {\em Computer Techniques for Electromagnetics},
  International Series of Monographs in Electrical Engineering, chapter~4,
  pages 159--264. Pergamon, Oxford, UK, 1973.

\bibitem{chang1974surface}
Yu~Chang and Roger~F Harrington.
\newblock A surface formulation for characteristic modes of material bodies.
\newblock Technical report, Syracuse University, Syracuse, NY, 1974.
\newblock Technical Report TR-74-7.

\bibitem{wu1977scattering-bor}
Te-Kao Wu and Leonard~L Tsai.
\newblock Scattering from arbitrarily-shaped lossy dielectric bodies of
  revolution.
\newblock {\em Radio Science}, 12(5):709--718, 1977.

\bibitem{arvas1986field}
E.~Arvas, S.~M. Rao, and T.~K. Sarkar.
\newblock {E}-field solution of {TM}-scattering from multiple perfectly
  conducting and lossy dielectric cylinders of arbitrary cross-section.
\newblock {\em IEE Proceedings H (Microwaves, Antennas and Propagation)},
  133(2):115--121, 1986.

\bibitem{wu2012fast}
Haijun Wu, Yijun Liu, and Weikang Jiang.
\newblock A fast multipole boundary element method for {3D} multi-domain
  acoustic scattering problems based on the {B}urton-{M}iller formulation.
\newblock {\em Engineering Analysis with Boundary Elements}, 36(5):779--788,
  2012.

\bibitem{harrington1989boundary}
Roger~F Harrington.
\newblock Boundary integral formulations for homogeneous material bodies.
\newblock {\em Journal of Electromagnetic Waves and Applications}, 3(1):1--15,
  1989.

\bibitem{rao1990field}
Sadasiva~M Rao and Donald~R Wilton.
\newblock {E}-field, {H}-field, and combined field solution for arbitrarily
  shaped three-dimensional dielectric bodies.
\newblock {\em Electromagnetics}, 10(4):407--421, 1990.

\bibitem{sikora2006diffuse}
Jan Sikora, Athanasios Zacharopoulos, Abdel Douiri, Martin Schweiger, Lior
  Horesh, Simon~R Arridge, and Jorge Ripoll.
\newblock Diffuse photon propagation in multilayered geometries.
\newblock {\em Physics in Medicine \& Biology}, 51(3):497, 2006.

\bibitem{juffer1991electric}
Andr{\'e}~H Juffer, Eugen~FF Botta, Bert~AM {van Keulen}, Auke {van der Ploeg},
  and Herman~JC Berendsen.
\newblock The electric potential of a macromolecule in a solvent: A fundamental
  approach.
\newblock {\em Journal of Computational Physics}, 97(1):144--171, 1991.

\bibitem{bardhan2009numerical}
Jaydeep~P Bardhan.
\newblock Numerical solution of boundary-integral equations for molecular
  electrostatics.
\newblock {\em The Journal of Chemical Physics}, 130(9):094102, 2009.

\bibitem{hiptmair2012multiple}
Ralf Hiptmair and C~Jerez-Hanckes.
\newblock Multiple traces boundary integral formulation for {H}elmholtz
  transmission problems.
\newblock {\em Advances in Computational Mathematics}, 37(1):39--91, 2012.

\bibitem{peng2012computations}
Zhen Peng, Kheng-Hwee Lim, and Jin-Fa Lee.
\newblock Computations of electromagnetic wave scattering from penetrable
  composite targets using a surface integral equation method with multiple
  traces.
\newblock {\em IEEE Transactions on Antennas and Propagation}, 61(1):256--270,
  2012.

\bibitem{claeys2015integral}
Xavier Claeys and Ralf Hiptmair.
\newblock Integral equations for acoustic scattering by partially impenetrable
  composite objects.
\newblock {\em Integral Equations and Operator Theory}, 81(2):151--189, 2015.

\bibitem{claeys2019introduction}
Xavier Claeys, Victorita Dolean, and Martin~J Gander.
\newblock An introduction to multi-trace formulations and associated domain
  decomposition solvers.
\newblock {\em Applied Numerical Mathematics}, 135:69--86, 2019.

\bibitem{marx1982single}
Egon Marx.
\newblock Single integral equation for wave scattering.
\newblock {\em Journal of Mathematical Physics}, 23(6):1057--1065, 1982.

\bibitem{marx1984integral}
Egon Marx.
\newblock Integral equation for scattering by a dielectric.
\newblock {\em IEEE Transactions on Antennas and Propagation}, 32(2):166--172,
  1984.

\bibitem{glisson1984integral}
A.~Glisson.
\newblock An integral equation for electromagnetic scattering from homogeneous
  dielectric bodies.
\newblock {\em IEEE Transactions on Antennas and propagation}, 32(2):173--175,
  1984.

\bibitem{mautz1989stable}
Joseph~R Mautz.
\newblock A stable integral equation for electromagnetic scattering from
  homogeneous dielectric bodies.
\newblock {\em IEEE Transactions on Antennas and Propagation},
  37(8):1070--1071, 1989.

\bibitem{kleinman1988single}
RE~Kleinman and PA~Martin.
\newblock On single integral equations for the transmission problem of
  acoustics.
\newblock {\em SIAM Journal on Applied Mathematics}, 48(2):307--325, 1988.

\bibitem{yla2011calderon}
Pasi Yl{\"a}-Oijala, Sami~P Kiminki, and Seppo J{\"a}rvenp{\"a}{\"a}.
\newblock {C}alder{\'o}n preconditioned surface integral equations for
  composite objects with junctions.
\newblock {\em IEEE Transactions on Antennas and Propagation}, 59(2):546--554,
  2011.

\bibitem{gossye2018calderon}
Michiel Gossye, Martijn Huynen, Dries {Vande Ginste}, Dani{\"e}l {De Zutter},
  and Hendrik Rogier.
\newblock A {C}alder{\'o}n preconditioner for high dielectric contrast media.
\newblock {\em IEEE Transactions on Antennas and Propagation}, 66(2):808--818,
  2018.

\bibitem{gossye2019electromagnetic}
Michiel Gossye, Dries {Vande Ginste}, and Hendrik Rogier.
\newblock Electromagnetic modeling of high magnetic contrast media using
  {C}alder{\'o}n preconditioning.
\newblock {\em Computers \& Mathematics with Applications}, 77(6):1626--1638,
  2019.

\bibitem{wout2022highcontrast}
Elwin van~'t Wout, Seyyed~R. Haqshenas, Pierre G{\'e}lat, Timo Betcke, and
  Nader Saffari.
\newblock Boundary integral formulations for acoustic modelling of
  high-contrast media.
\newblock {\em Computers \& Mathematics with Applications}, 105:136--149, 2022.

\bibitem{rapun2006indirect}
Mar{\'\i}a-Luisa Rap{\'u}n and Francisco-Javier Sayas.
\newblock Indirect methods with {B}rakhage-{W}erner potentials for {H}elmholtz
  transmission problems.
\newblock In {\em Numerical Mathematics and Advanced Applications}, pages
  1146--1154. Springer, 2006.

\bibitem{boubendir2015integral}
Yassine Boubendir, Oscar Bruno, David Levadoux, and Catalin Turc.
\newblock Integral equations requiring small numbers of {K}rylov-subspace
  iterations for two-dimensional smooth penetrable scattering problems.
\newblock {\em Applied Numerical Mathematics}, 95:82--98, 2015.

\bibitem{claeys2013multi}
Xavier Claeys and Ralf Hiptmair.
\newblock Multi-trace boundary integral formulation for acoustic scattering by
  composite structures.
\newblock {\em Communications on Pure and Applied Mathematics},
  66(8):1163--1201, 2013.

\bibitem{laliena2009symmetric}
Antonio~R Laliena, M-L Rap{\'u}n, and F-J Sayas.
\newblock Symmetric boundary integral formulations for {H}elmholtz transmission
  problems.
\newblock {\em Applied Numerical Mathematics}, 59(11):2814--2823, 2009.

\bibitem{langer2003boundary}
Ulrich Langer and Olaf Steinbach.
\newblock Boundary element tearing and interconnecting methods.
\newblock {\em Computing}, 71(3):205--228, 2003.

\bibitem{wu2008multi}
Ting-Wen Wu.
\newblock Multi--domain boundary element method in acoustics.
\newblock In {\em Computational Acoustics of Noise Propagation in Fluids-Finite
  and Boundary Element Methods}, pages 367--386. Springer, 2008.

\bibitem{peng2012nonconformal}
Zhen Peng, Kheng-Hwee Lim, and Jin-Fa Lee.
\newblock Nonconformal domain decomposition methods for solving large
  multiscale electromagnetic scattering problems.
\newblock {\em Proceedings of the IEEE}, 101(2):298--319, 2012.

\bibitem{xia2016enhanced}
Tian Xia, Hui Gan, Michael Wei, Weng~Cho Chew, Henning Braunisch, Zhiguo Qian,
  Kemal Ayg{\"u}n, and Alaeddin Aydiner.
\newblock An enhanced augmented electric-field integral equation formulation
  for dielectric objects.
\newblock {\em IEEE Transactions on Antennas and Propagation},
  64(6):2339--2347, 2016.

\bibitem{chew2014vector}
Weng~Cho Chew.
\newblock Vector potential electromagnetics with generalized gauge for
  inhomogeneous media: Formulation.
\newblock {\em Progress In Electromagnetics Research}, 149:69--84, 2014.

\bibitem{bremer2015high}
James Bremer, Adrianna Gillman, and Per-Gunnar Martinsson.
\newblock A high-order accurate accelerated direct solver for acoustic
  scattering from surfaces.
\newblock {\em BIT Numerical Mathematics}, 55(2):367--397, 2015.

\bibitem{marburg2003performance}
Steffen Marburg and Stefan Schneider.
\newblock Performance of iterative solvers for acoustic problems. {P}art {I}.
  {S}olvers and effect of diagonal preconditioning.
\newblock {\em Engineering Analysis with Boundary Elements}, 27(7):727--750,
  2003.

\bibitem{saad1986gmres}
Youcef Saad and Martin~H Schultz.
\newblock Gmres: A generalized minimal residual algorithm for solving
  nonsymmetric linear systems.
\newblock {\em SIAM Journal on Scientific and Statistical Computing},
  7(3):856--869, 1986.

\bibitem{schneider2003performance}
Stefan Schneider and Steffen Marburg.
\newblock Performance of iterative solvers for acoustic problems. {P}art {II}.
  {A}cceleration by {ILU}-type preconditioner.
\newblock {\em Engineering Analysis with Boundary Elements}, 27(7):751--757,
  2003.

\bibitem{sakuma2008fast}
Tetsuya Sakuma, Stefan Schneider, and Yosuke Yasuda.
\newblock Fast solution methods.
\newblock In {\em Computational Acoustics of Noise Propagation in Fluids-Finite
  and Boundary Element Methods}, chapter~12, pages 333--366. Springer, 2008.

\bibitem{darbas2013combining}
Marion Darbas, Eric Darrigrand, and Yvon Lafranche.
\newblock Combining analytic preconditioner and fast multipole method for the
  3-{D} {H}elmholtz equation.
\newblock {\em Journal of Computational Physics}, 236:289--316, 2013.

\bibitem{antoine2008integral}
Xavier Antoine and Yassine Boubendir.
\newblock An integral preconditioner for solving the two-dimensional scattering
  transmission problem using integral equations.
\newblock {\em International Journal of Computer Mathematics},
  85(10):1473--1490, 2008.

\bibitem{hiptmair2006operator}
Ralf Hiptmair.
\newblock Operator preconditioning.
\newblock {\em Computers \& Mathematics with Applications}, 52(5):699--706,
  2006.

\bibitem{kirby2010functional}
Robert~C Kirby.
\newblock From functional analysis to iterative methods.
\newblock {\em SIAM Review}, 52(2):269--293, 2010.

\bibitem{betcke2020product}
Timo Betcke, Matthew~W Scroggs, and Wojciech {\'S}migaj.
\newblock Product algebras for {G}alerkin discretisations of boundary integral
  operators and their applications.
\newblock {\em ACM Transactions on Mathematical Software (TOMS)}, 46(1):1--22,
  2020.

\bibitem{yan2010comparative}
Su~Yan, Jian-Ming Jin, and Zaiping Nie.
\newblock A comparative study of {C}alder{\'o}n preconditioners for {PMCHWT}
  equations.
\newblock {\em IEEE Transactions on Antennas and Propagation},
  58(7):2375--2383, 2010.

\bibitem{steinbach1998construction}
Olaf Steinbach and Wolfgang~L Wendland.
\newblock The construction of some efficient preconditioners in the boundary
  element method.
\newblock {\em Advances in Computational Mathematics}, 9(1-2):191--216, 1998.

\bibitem{moore1988theory}
Thomas~G Moore, Jeffery~G Blaschak, Allen Taflove, and Gregory~A Kriegsmann.
\newblock Theory and application of radiation boundary operators.
\newblock {\em IEEE Transactions on Antennas and Propagation},
  36(12):1797--1812, 1988.

\bibitem{antoine2008advances}
Xavier Antoine.
\newblock Advances in the on-surface radiation condition method: Theory,
  numerics and applications.
\newblock {\em Computational Methods for Acoustics Problems}, pages 169--194,
  2008.

\bibitem{antoine2005alternative}
Xavier Antoine and Marion Darbas.
\newblock Alternative integral equations for the iterative solution of acoustic
  scattering problems.
\newblock {\em The Quarterly Journal of Mechanics and Applied Mathematics},
  58(1):107--128, 2005.

\bibitem{haqshenas2020fast}
S.~R. Haqshenas, P.~G{\'e}lat, E.~{van 't Wout}, T.~Betcke, and N.~Saffari.
\newblock A fast full-wave solver for calculating ultrasound propagation in the
  body.
\newblock {\em Ultrasonics}, page 106240, 2020.

\bibitem{wout2022pmchwt}
Elwin van~'t Wout, Seyyed~R. Haqshenas, Pierre G{\'e}lat, Timo Betcke, and
  Nader Saffari.
\newblock Frequency-robust preconditioning of boundary integral equations for
  acoustic transmission.
\newblock {\em Journal of Computational Physics}, 462:111229, 2022.

\bibitem{wout2015fast}
Elwin {van 't Wout}, Pierre G{\'e}lat, Timo Betcke, and Simon Arridge.
\newblock A fast boundary element method for the scattering analysis of
  high-intensity focused ultrasound.
\newblock {\em The Journal of the Acoustical Society of America},
  138(5):2726--2737, 2015.

\bibitem{betcke2017computationally}
Timo Betcke, Elwin {van 't Wout}, and Pierre G{\'e}lat.
\newblock Computationally efficient boundary element methods for high-frequency
  {H}elmholtz problems in unbounded domains.
\newblock In Domenico Lahaye, Jok Tang, and Kees Vuik, editors, {\em Modern
  Solvers for {H}elmholtz Problems}, pages 215--243. Birkh\"auser, Cham, 2017.

\bibitem{wout2021proximity}
Elwin {van 't Wout} and Christopher Feuillade.
\newblock Proximity resonances of water-entrained air bubbles near acoustically
  reflecting boundaries.
\newblock {\em The Journal of the Acoustical Society of America},
  149(4):2477--–2491, 2021.

\bibitem{hornikx2015platform}
Maarten Hornikx, Manfred Kaltenbacher, and Steffen Marburg.
\newblock A platform for benchmark cases in computational acoustics.
\newblock {\em Acta Acustica United With Acustica}, 101(4):811--820, 2015.

\bibitem{hamilton1998nonlinear}
Mark~F Hamilton and David~T Blackstock.
\newblock {\em Nonlinear acoustics}, volume 237.
\newblock Academic Press, 1998.

\bibitem{duck1990physical}
FA~Duck.
\newblock {\em Physical properties of tissue: a comprehensive reference book}.
\newblock Academic Press, London, UK, 1990.

\bibitem{itis2018}
{IT'IS Foundation}.
\newblock Tissue properties database.
\newblock Online database, 2018.

\bibitem{geuzaine2009gmsh}
Christophe Geuzaine and Jean-Fran{\c{c}}ois Remacle.
\newblock Gmsh: A {3-D} finite element mesh generator with built-in pre-and
  post-processing facilities.
\newblock {\em International Journal for Numerical Methods in Engineering},
  79(11):1309--1331, 2009.

\bibitem{zhang2020dual}
Jianming Zhang, Weicheng Lin, Xiaomin Shu, and Yudong Zhong.
\newblock A dual interpolation boundary face method for exterior acoustic
  problems based on the {B}urton-{M}iller formulation.
\newblock {\em Engineering Analysis with Boundary Elements}, 113:219--231,
  2020.

\bibitem{marburg2002six}
Steffen Marburg.
\newblock Six boundary elements per wavelength: Is that enough?
\newblock {\em Journal of Computational Acoustics}, 10(01):25--51, 2002.

\bibitem{scroggs2017software}
Matthew~W Scroggs, Timo Betcke, Erik Burman, Wojciech {\'S}migaj, and Elwin
  {van 't Wout}.
\newblock Software frameworks for integral equations in electromagnetic
  scattering based on {C}alder{\'o}n identities.
\newblock {\em Computers \& Mathematics with Applications}, 74(11):2897--2914,
  2017.

\bibitem{scipy}
Pauli Virtanen, Ralf Gommers, Travis~E. Oliphant, Matt Haberland, Tyler Reddy,
  David Cournapeau, Evgeni Burovski, Pearu Peterson, Warren Weckesser, Jonathan
  Bright, St{\'e}fan~J. {van der Walt}, Matthew Brett, Joshua Wilson, K.~Jarrod
  Millman, Nikolay Mayorov, Andrew R.~J. Nelson, Eric Jones, Robert Kern, Eric
  Larson, C~J Carey, {\.I}lhan Polat, Yu~Feng, Eric~W. Moore, Jake
  {VanderPlas}, Denis Laxalde, Josef Perktold, Robert Cimrman, Ian Henriksen,
  E.~A. Quintero, Charles~R. Harris, Anne~M. Archibald, Ant{\^o}nio~H. Ribeiro,
  Fabian Pedregosa, Paul {van Mulbregt}, and {SciPy 1.0 Contributors}.
\newblock {{SciPy} 1.0: Fundamental Algorithms for Scientific Computing in
  {P}ython}.
\newblock {\em Nature Methods}, 17:261--272, 2020.

\bibitem{hunter2007matplotlib}
J.~D. Hunter.
\newblock Matplotlib: A {2D} graphics environment.
\newblock {\em Computing in Science \& Engineering}, 9(3):90--95, 2007.

\bibitem{waskom2020seaborn}
Michael Waskom.
\newblock mwaskom/seaborn.
\newblock Zenodo, September 2020.

\bibitem{molina2018iterative}
Jorge Molina-Moya, Alejandro~Enrique Mart{\'\i}nez-Castro, and Pablo Ortiz.
\newblock An iterative parallel solver in {GPU} applied to frequency domain
  linear water wave problems by the boundary element method.
\newblock {\em Frontiers in Built Environment}, 4:69, 2018.

\bibitem{yokota2012tuned}
Rio Yokota and Lorena~A Barba.
\newblock A tuned and scalable fast multipole method as a preeminent algorithm
  for exascale systems.
\newblock {\em The International Journal of High Performance Computing
  Applications}, 26(4):337--346, 2012.

\bibitem{golub2013matrix}
G.H. Golub and C.F. {Van Loan}.
\newblock {\em Matrix Computations}.
\newblock Johns Hopkins University Press, Baltimore, MD, 2013.

\bibitem{heike2017letter}
Hofmann Heike, H~Wickham, and K~Kafadar.
\newblock Letter-value plots: Boxplots for large data.
\newblock {\em Journal of Computational and Graphical Statistics}, 26:469--477,
  2017.

\end{thebibliography}

\end{document}